\documentclass{article}%
\usepackage{amssymb}
\usepackage{amsmath}
\usepackage{amsfonts}
\usepackage{graphicx}%
\providecommand{\U}[1]{\protect\rule{.1in}{.1in}}

\newenvironment{proof}[1][Proof]{\noindent\textbf{#1.} }{\ \rule{0.5em}{0.5em}}
\begin{document}

\title{ }

\begin{center}
\bigskip{\LARGE A solution of the maximality problem for one-parameter
dynamical systems }

\bigskip

Costel Peligrad

Department of Mathematical Sciences, University of Cincinnati

4508 French Hall West, Cincinnati, Ohio 45221-0025

Email address: costel.peligrad@uc.edu
\end{center}

\bigskip

\textbf{Abstract \ }We prove a maximality theorem for one-parameter dynamical
systems that include W*- , C*- and multiplier one-parameter dynamical systems.
Our main result is new even for one-parameter actions on commutative
multiplier algebras including the algebra $C_{b}(%
\mathbb{R}
)$ of bounded continuous functions on $%
\mathbb{R}
$\ acted upon by translations. The methods we develop and use in our
characterization of maximality include harmonic analysis, topological vector
spaces and operator algebra techniques.

\section{ Introduction}

The current paper grew out of an effort to find a solution to the general
problem of maximality of subalgebras of analytic elements associated to
various dynamical systems: W*-dynamical systems, C*-dynamical systems and, as
it will follow from the present work, multiplier dynamical systems as defined
below (Section 2). The study of maximality of analytic subalgebras associated
with C*- or W*-dynamical systems has a history of over six decades. It started
with Wermer's maximality theorem [27]. Motivated by some earlier problems
about approxmation of continuous functions, Wermer showed that if
$X=C(\mathbf{T})$\ is the C*-algebra of all continuous functions on
$\mathbf{T}=\left\{  t\in\mathbf{%
\mathbb{C}
:}\left\vert t\right\vert =1\right\}  ,$\ then the norm closed subalgebra
$A$\ of all functions $f\in C(\mathbf{T})$\ that have an analytic extension to
the unit disc $\mathbf{D}=\left\{  s\in\mathbf{%
\mathbb{C}
:}\left\vert s\right\vert <1\right\}  $\ is a maximal norm-closed subalgebra
of $X.$ Hoffman and Singer [9] and Simon [23] obtained some generalizations of
Wermer's theorem for the case of compact groups with archimedean-linearly
ordered Pontryagin duals.

A significant and far reaching generalization of Wermer's maximality theorem
was obtained by Forelli [7] using the seminal concepts and results from his
paper [6]. In [7], Forelli considered a minimal action of $%
\mathbb{R}
$ on a locally compact Hausdorff space $S,$\ that is a homomorphism $\tau$\ of
$%
\mathbb{R}
$\ into the group of homeomorphisms of $S$\ onto itself such that the mapping
$t\rightarrow\tau_{t}(s)$\ is continuous for every $s\in S$\ and every orbit
$\left\{  \tau_{t}(s):t\in%
\mathbb{R}
\right\}  ,s\in S$\ is dense in $S.$ If $f\in C_{0}(S)$\ and $t\in%
\mathbb{R}
,$ \ denote $\alpha_{t}(f)=f\circ\tau_{t}.$\ Then $\alpha$\ is a homeomorphism
of $%
\mathbb{R}
$\ into the group of automorphisms of the C*-algebra $X=C_{0}(S)$\ such that
the mapping $t\rightarrow\alpha_{t}(f)$\ is continuous from $%
\mathbb{R}
$ to $C_{0}(S)$ endowed with the uniform norm topology\ for every $f\in
C_{0}(S).$ Forelli proved that the subalgebra $X^{\alpha}([0,\infty))$\ of
$X=C_{0}(S)$\ consisting of all functions $f\in C_{0}(S)$\ such that the
$\alpha-$spectrum of $f,$ $sp_{\alpha}(f)$,\ as defined in the next
section,\ is contained in $[0,\infty),$ is a maximal norm closed subalgebra of
$X.$ An equivalent description of the algebra $X^{\alpha}([0,\infty))$\ for
$X=C_{0}(S)$ and $\alpha$ as defined above\ is the following: $f\in X^{\alpha
}([0,\infty))$\ if and only if the mapping $t\rightarrow\alpha_{t}%
(f)(s)=(f\circ\tau_{t})(s)$\ has a bounded analytic extension to the upper
half plane for every $s\in S.$ In what follows a system $(X,G,\alpha
)$\ consisting of a C* algebra $X,$ a locally compact group $G$\ and a
homeomorphism $\alpha$\ from $G$\ into the group $Aut(X)$\ of all
automorphisms of $X$\ such that the mapping $t\rightarrow\alpha_{t}(x)$\ is
continuous from $G$ to $X$ with the norm topology\ for every $x\in X,$ will be
called a C*-dynamical system. When $G=%
\mathbb{R}
$\ the system will be called a one-parameter C*-dynamical system. If $X$\ does
not contain any non trivial norm-closed $\alpha$-invariant ideal, the system
will be called $\alpha$-simple. This is the case with Forelli system
$(C_{0}(S),%
\mathbb{R}
,\alpha)$ if $\tau$ is a minimal flow as defined above.

Another direction of study of maximal subalgebras is the following: If $(X,%
\mathbb{R}
,\alpha)$ is a W*-dynamical system consisting of a von Neumann algebra $X$ and
an action of $%
\mathbb{R}
$ on $X$\ such that the mapping $t\rightarrow\alpha_{t}(x)$\ is continuous
from $%
\mathbb{R}
$\ to $X$\ endowed with the $w$*-topology, for every $x\in X,$ when is
$X^{\alpha}([0,\infty))$\ a maximal $w$*-closed subalgebra of $X$? In [21]
Sarason noticed that Hoffman and Singer [10] basically proved the first result
in this direction: $H^{\infty}(\mathbf{T})$ is a maximal $w$*-closed
subalgebra of $L^{\infty}(\mathbf{T}),$ where $H^{\infty}(\mathbf{T})\subset
L^{\infty}(\mathbf{T})$ is the subalgebra of all $f\in L^{\infty}(\mathbf{T})$
that have an analytic extension to the unit disk $\mathbf{D.}$ Further in [14]
Muhly showed that if $(X,%
\mathbb{R}
,\alpha)$ is an ergodic W*-dynamical system with $X$\ an abelian von Neumann
algebra, then the analytic subalgebra, $X^{\alpha}([0,\infty))$ as defined
below is a maximal w*-closed subalgebra of $X.$

\ \ Some extensions of Wermer's result to the more general case of C*- and
W*-dynamical systems were obtained in [13] for W*-crossed products and in [18]
for C*-crossed products. In [24] Solel has found necessary and sufficient
conditions for maximality of subalgebras of analytic elements associated with
periodic $\sigma$-finite W*-dynamical systems and, further, in [25] he has
considered the maximality for the particular case of an $\alpha$-simple
one-parameter $\sigma$-finite W*-dynamical system, i.e. the case when $X$\ has
no non-trivial $w^{\ast}$-closed $\alpha$-invariant ideals (or equivalently
the center of $X\ $ contains no non-trivial $\alpha$-invariant projections).

In [19] we found a necessary and sufficient condition, that we called
Condition (\textbf{S}) (\textbf{S} from the word spectrum), for maximality of
the analytic subalgebra associated with a periodic C*-dynamical system. Soon
after, Kishimoto [12] considered the case of one-parameter C*-dynamical
systems $(X,%
\mathbb{R}
,\alpha).$ He proved that if the crossed product C*-algebra of the system is a
simple C*-algebra, then $X^{\alpha}([0,\infty))$ is a maximal norm-closed
subalgebra of $X$.

In this paper we will find a spectral characterization of maximality for
general one-parameter dynamical systems, not necessarily separable or $\sigma
$-finite (Condition 4.3.). This characterization is equivalent to our
Condition (\textbf{S)} for periodic C*-dynamical systems [19] and with Solel
characterization of maximality for $\sigma$-finite W*-dynamical systems [24].
It characterizes also maximality for one-parameter multiplier dynamical
systems (as defined in Remark 2.2.). The methods we develop and use include
harmonic analysis, topological vector spaces and operator algebra techniques.
Our main result\ (Theorem 4.18.) contains and improves on all of the above
mentioned results. In particular it contains the following result that was not
considered so far even for commutative dynamical systems: Let $C_{b}(%
\mathbb{R}
)$ be the C*-algebra of bounded continuous functions\ on $%
\mathbb{R}
$ and $\alpha_{t}(f)(s)=f(s-t)$ be the translation of $f$\ by $t\in%
\mathbb{R}
.$ Then $C_{b}(%
\mathbb{R}
)$ is the multiplier C*-algebra of $C_{0}(%
\mathbb{R}
)$ and for each $f\in C_{b}(%
\mathbb{R}
)$ the mapping $t\rightarrow\alpha_{t}(f)$ is continuous from $%
\mathbb{R}
$ to $C_{b}(%
\mathbb{R}
)$ endowed with the strict topology. A consequence of our main result is that
$H^{\infty}(%
\mathbb{R}
)\cap C_{b}(%
\mathbb{R}
)$ is a maximal strictly closed subalgebra of $C_{b}(%
\mathbb{R}
).$

This paper is organized as follows. In Section 2 we discuss the concept of
dual pairs, $(X,\mathcal{F})$ of Banach spaces ([1], [5], [28]) and
representations of locally compact groups on a Banach space $X$ endowed with
the $\mathcal{F}$-topology$.$\ We will define the spectral subspaces
associated with these representations and define the $\mathcal{F}$-dynamical
systems. We also define the concept of strong $\mathcal{F}$-Connes spectrum of
such systems and\ we show that our notion of strong $\mathcal{F}$-Connes
spectrum contains both the original concept of Connes spectrum for
W*-dynamical systems and the strong Connes spectrum defined by Kishimoto for
C*-dynamical systems. The strong $\mathcal{F}$-Connes spectrum is a new
concept for multiplier dynamical systems. In Section 3 we obtain some results
about periodic $\mathcal{F}$-dynamical systems that will be used in Section 4.
These results extend to the case of $\mathcal{F}$-dynamical systems results by
Connes [3] for W*-dynamical systems and Olesen, Pedersen and Stormer [16] for
C*-dynamical systems. In Section 4 we prove our main result about maximality
of subalgebras of analytic elements of an $\mathcal{F}$-dynamical system. I am
grateful to Laszlo Zsido for several useful discussions and suggestions during
the completion of this work. I also thank Akitaka Kishimoto for enlightening
discussions about his work on the subject.

\section{\bigskip Notations and preliminary results}

In this section we will set up the framework for fututure discussions and will
establish the notations that will be used in the rest of the paper. We will
also state and prove some basic preliminary results.

\textbf{ }

\textbf{2.1. Definition. }\textit{A dual pair of Banach spaces is, by
definition a pair (}$X,\mathcal{F})$\textit{ of Banach spaces together with a
bilinear functional}%
\[
(x,\varphi)\in X\times\mathcal{F}\rightarrow\left\langle x,\varphi
\right\rangle =\varphi(x)\in%
\mathbb{C}
\]
\textit{such that \newline1) }$\left\Vert x\right\Vert =\sup_{\varphi
\in\mathcal{F},\left\Vert \varphi\right\Vert \leq1}\left\vert \left\langle
x,\varphi\right\rangle \right\vert ,x\in X$\textit{\newline2) }$\left\Vert
\varphi\right\Vert =\sup_{x\in X,\left\Vert x\right\Vert \leq1}\left\vert
\left\langle x,\varphi\right\rangle \right\vert ,\varphi\in F$\textit{\newline%
3) The convex hull of every relatively }$\mathcal{F}$\textit{-compact subset
of }$X$\textit{\ is relatively }$\mathcal{F}$\textit{-compact\newline4) The
convex hull of every relatively }$X$\textit{-compact subset of }%
$F$\textit{\ is relatively }$X$\textit{-compact.}\newline\textit{If }%
$X$\textit{\ is a C*-algebra we will assume in addition that the following
conditions hold}\newline\textit{5) The involution of }$X$\textit{\ is
}$\mathcal{F}$\textit{-continuous and the multiplication in }$X$\textit{\ is
separately }$F$\textit{-continuous.}

\bigskip

Clearly, if \textit{(}$X,\mathcal{F})$ is a dual pair of Banach spaces as in
the above definition, then $X$\ is naturally embedded in $\mathcal{F}^{\ast}$
and $\mathcal{F}$\ is naturally embedded in $X^{\ast}.$

\bigskip

\textbf{2.2. Remark} \textit{i) If }$X$\textit{ is a C*-algebra, and
}$\mathcal{F}=X^{\ast}$\textit{ is its dual then the pair }$(X,X^{\ast}%
)$\textit{ satisfies conditions 1)-5)\newline ii) If }$X$\textit{ is a
W*-algebra and }$\mathcal{F}=X_{\ast}$\textit{ is its predual, then the pair
}$(X,X_{\ast})$\textit{ satisfies conditions 1)-5)\newline iii) If }%
$X=M(Y)$\textit{\ is the multiplier algebra of a C*-algebra }$Y$\textit{\ and
}$\mathcal{F}=Y^{\ast}$\textit{\ is the dual of }$Y,$\textit{\ then the pair
}$(X,Y^{\ast})$\textit{ satisfies conditions 1)-5). In addition, in this case,
the duality is compatible with the strict topology on }$X=M(Y).$

\bigskip

\begin{proof}
Parts i) and ii) are discussed in [1]. Part iii) follows from [4] and [26].
\end{proof}

\bigskip

Notice that in all three cases i), ii) and iii) of the above Remark 2.2., the
dual $\mathcal{F}^{\ast}$\ of $\mathcal{F}$\ is a von Neumann algebra and
$X$\ is naturally embedded in $\mathcal{F}^{\ast}.$

\bigskip

\textbf{2.3. Remark }\textit{i) Let (}$X,\mathcal{F})$\textit{ be a dual pair
of Banach spaces with }$\mathcal{F}=X^{\ast}$\textit{. If }$Y\subset
X$\textit{ is a norm closed subspace of }$X$\textit{, then }$Y$\textit{\ is
}$\mathcal{F}-$\textit{closed in }$X$\textit{.\newline ii) Let (}%
$X,\mathcal{F})$\textit{ be a dual pair of Banach spaces. If }$Y\subset
X$\textit{ is an }$\mathcal{F}$\textit{-closed subspace, then }$Y$\textit{\ is
norm closed.}

\bigskip

\begin{proof}
Part i) follows from the Hahn Banach separation theorem. Part ii) is trivial
since the $\mathcal{F}$-topology is weaker than the norm topology.
\end{proof}

\bigskip

\textbf{2.4. Notation} \textit{Let (}$X,\mathcal{F})$\textit{ be a dual pair
of Banach spaces. }I\textit{f }$Y,Z$\textit{\ are subsets of }$X$%
\textit{\ denote:}

\textit{a) }$lin\left\{  Y\right\}  $\textit{\ is the linear span of }%
$Y$\textit{.}

\textit{b) If }$X\ $\textit{is a Banach space with involution, in particular a
C*-algebra, denote }$Y^{\ast}=\left\{  y^{\ast}:y\in Y\right\}  .$

\textit{c) If }$X$\textit{ is a C*-algebra, denote }$YZ=lin\left\{  yz:y\in
Y,z\in Z\right\}  .$

\textit{d) }$\overline{Y}^{\sigma}=\mathcal{F}-$\textit{closure of }$Y$ in
$X.$

\textit{e) }$\overline{Y}^{\left\Vert {}\right\Vert }=$\textit{ norm closure
of }$Y.$

\textit{f)} $\overline{Y}^{w}=$ \textit{the }$w^{\ast}$\textit{-closure of
}$Y$\textit{\ in} $\mathcal{F}^{\ast}.$

\textit{Suppose that (}$X,\mathcal{F})$\textit{ is a dual pair of Banach
spaces with }$X$\textit{\ a C*-algebra, so Definition 2.1. 1)-5) is
satisfied.}

\textit{g) If }$\varphi\in\mathcal{F},$ \textit{denote by} $\varphi^{\ast}\in
X^{\ast}$\textit{ the adjoint} of $\varphi,$ \textit{defined by}
$\varphi^{\ast}(x)=\overline{\varphi(x^{\ast})},x\in X,$ \textit{where
}$\overline{\varphi(x^{\ast})}$\textit{ is the complex conjugate of }%
$\varphi(x^{\ast}).$

\textit{h) If }$\varphi\in F$\textit{ and }$y\in X$\textit{ denote by }%
$L_{y}\varphi\in X^{\ast}$\textit{ the functional defined by (}$L_{y}%
\varphi)(x)=\varphi(yx),x\in X.$

\bigskip

Suppose that $(X,\mathcal{F})$\ is a dual pair of Banach spaces satisfying
conditions 1)-4) of the Definition 2.1. and $G$ a locally compact abelian group.

\bigskip

\textbf{2.5. Definition [1]} \textit{A representation of }$G$\textit{ on }%
$X$\textit{\ is a homomorphism of }$G$\textit{\ into the group of all
invertible bounded }$\mathcal{F}$\textit{-continuous linear operators on }%
$X$\textit{ such that }$\sup_{t\in G}\left\Vert \alpha_{t}\right\Vert <\infty
$\textit{ and for every }$x\in X,$\textit{ and }$\varphi\in\mathcal{F}%
,$\textit{\ the map }$t\rightarrow\varphi(\alpha_{t}(x))$\textit{ is
continuous. The triple }$(\mathcal{F},G,\alpha^{\ast})$\textit{\ is called the
dual representation of }$G$\textit{\ on }$\mathcal{F}.$

\bigskip

If $X$ is a C*-algebra and $\alpha_{t},t\in G$ are $\mathcal{F}$-continuous
automorphisms of $X,$\ the triple ($X,G,\alpha$) will be called an
$\mathcal{F}$-dynamical system. If $G=%
\mathbb{R}
,$\ the system ($X,%
\mathbb{R}
,\alpha)$ will be called a one-parameter $\mathcal{F}$-dynamical system.

\bigskip

\textbf{2.6. Notation.} \textit{If }$X$\textit{\ is a C*-algebra, we will also
use the following additional notations}

\textit{g) }$\mathcal{H}_{\sigma}^{\alpha}(X)$\textit{ the set of all non-zero
globally }$\alpha-$\textit{invariant }$\mathcal{F}$\textit{-closed hereditary
C*-subalgebras of }$X.$

\textit{h) }$\mathcal{H}^{\alpha}(X)$\textit{ the set of all non-zero globally
}$\alpha-$\textit{invariant hereditary C*-subalgebras of }$X.$\textit{\newline%
\ \ \ \ \ i) }$\mathcal{H}_{\sigma}(X)$\textit{ the set of all non-zero
}$\mathcal{F}$\textit{-closed hereditary C*-subalgebras of }$X.$%
\textit{\newline\ \ \ \ \ j) }$\mathcal{H}(X)$\textit{ the set of all non-zero
norm-closed hereditary C*-subalgebras of }$X.$

\bigskip

\textbf{2.7. Remark} \textit{i) If }$\mathcal{F}=X^{\ast}$\textit{\ then the
above concept of }$\mathcal{F}$\textit{-dynamical system coincides with the
concept of C*-dynamical system. Indeed in this case, every automorphism of
}$X$\textit{ is }$\mathcal{F}$\textit{-continuous and the continuity of
}$g\rightarrow\varphi(\alpha_{g}(a))$\textit{ for every }$\varphi\in X^{\ast}%
$\textit{\ and }$a\in X$\textit{\ is equivalent with the continuity of
}$g\rightarrow\alpha_{g}(a)$\textit{\ for every }$a\in X,$\textit{\ in the
norm of }$X$ \textit{[8, p.306].\newline ii) If }$X$\textit{ is a W*-algebra
and }$\mathcal{F}=X_{\ast},$\textit{ then the concept of }$\mathcal{F}%
$\textit{-dynamical system coincides with the usual notion of W*-dynamical
system.}

\bigskip

\textbf{2.8. Lemma} \textit{Let }$(X,\mathcal{F})$\textit{ be a dual pair of
Banach spaces with }$X$\textit{\ a C*-algebra. If }$Y$\textit{ is a hereditary
C*-subalgebra of }$X,$\textit{ then }$\overline{Y}^{\sigma}$\textit{ is}%
$.$\textit{an }$\mathcal{F}-$\textit{closed hereditary C*-subalgebra of }$X.$

\bigskip

\begin{proof}
Since the involution is $\mathcal{F}$-continuous and the multiplication is
separately $\mathcal{F}$-continuous (Definition 2.1. 5)), it follows that
$\overline{Y}^{\sigma}$\ is a C*-algebra. We have to prove that $\overline
{Y}^{\sigma}$ is a hereditary C*-subalgebra of $X.$\ Let $y\in\overline
{Y}^{\sigma},y\geqslant0$ and $x\in X$ be such that $0\leq x\leq y.$ Since
$\overline{Y}^{\sigma}$ is a C*-subalgebra of $X,$\ $y^{\frac{1}{4}}%
\in\overline{Y}^{\sigma}.$ By [17, Proposition 1.4.5.] there exists $u\in X$
such that $x^{\frac{1}{2}}=u%
\operatorname{y}%
^{\frac{1}{4}}=%
\operatorname{y}%
^{\frac{1}{4}}u^{\ast}.$ Thus $x=%
\operatorname{y}%
^{\frac{1}{4}}%
\operatorname{u}%
^{\ast}u%
\operatorname{y}%
^{\frac{1}{4}}.$ Applying again Definition 2.1. 5), it follows that
$x\in\overline{YXY}^{\sigma}.$ But, as $Y$\ is a hereditary C*-subalgebra of
$X,$\ we have $YXY\subset Y.$\ Therefore $x\in\overline{Y}^{\sigma}.$
\end{proof}

\bigskip

\textbf{2.9. Lemma }\textit{Let }$(X,\mathcal{F})$\textit{ be a dual pair of
Banach spaces with }$X$\textit{\ a C*-algebra. If }$A\subset X$\textit{ is a
subset then, }$\overline{A^{\ast}XA}^{\sigma}$\textit{ is an }$\mathcal{F}%
$\textit{-closed hereditary subalgebra of }$X$\textit{\ which coincides with
}$\overline{A^{\ast}AXA^{\ast}A}^{\sigma}.$\textit{ Moreover, if }$Y$\textit{
is an }$\mathcal{F}$\textit{-closed hereditary subalgebra of }$X$%
\textit{\ such that }$A^{\ast}A\subset Y,$\textit{ then }$\overline{A^{\ast
}XA}^{\sigma}\subset Y.$

\textit{ }

\begin{proof}
If $p$ is the support projection of $A$\ in $X^{\ast\ast}$\ then, clearly,
$p$\ is the range projection of $A^{\ast}$\ and therefore the support and
range projection of $A^{\ast}A.$ Hence $\overline{A^{\ast}XA}^{\left\Vert
{}\right\Vert }=\overline{A^{\ast}AXA^{\ast}A}^{\left\Vert {}\right\Vert },$
so the first part of the Lemma is proven. To prove the second part, notice
that if $A^{\ast}A\subset Y,$ then $A^{\ast}AXA^{\ast}A\subset YXY\subset
\overline{YXY\ }^{\sigma}=Y,$ so $\overline{A^{\ast}XA}^{\sigma}%
=\overline{A^{\ast}AXA^{\ast}A}^{\sigma}\subset Y.$
\end{proof}

\bigskip

Let $(X,\mathcal{F})$ be a dual pair of Banach spaces, $\alpha\ $a
representation of the locally compact abelian group $G\ $on $X$ as defined
above.Then, for every $x\in X,$\ and $f\in L^{1}(G)$\ there exists a unique
$y=\alpha_{f}(x)\in X$\ also denoted by $\alpha_{f}(x)=\int_{G}f(t)\alpha
_{t}(x)dt,$\ such that
\[
\varphi(y)=\int_{G}f(t)\varphi(\alpha_{t}(x))dt.
\]
for every $\varphi\in\mathcal{F}$ and the operator $x\rightarrow\alpha
_{f}(x)\ $is $\mathcal{F}-$continuous [1, Proposition 1.4.]$.$ If $x\in
X,$\ define the Arveson\ spectrum of $x$%
\[
sp_{\alpha}(x)=\left\{  \gamma\in\widehat{G}:\widehat{f}(\gamma)=0\text{ for
all }f\in L^{1}(G)\text{\ with }\alpha_{f}(x)=0\right\}  .
\]
where $\widehat{f}$\ is the Fourier transform of $f.$\ Then, it is clear that%
\[
sp_{\alpha}(x^{\ast})=\left\{  -\gamma:\gamma\in sp_{\alpha}(x)\right\}
\]
The spectral subspaces of the system ($X,G,\alpha$)\ are defined as follows
([1], [17]).

If $F$\ is a closed subset of $\widehat{G},$denote
\[
X^{\alpha}(F)=\left\{  x\in X:sp_{\alpha}(x)\subset F\right\}  .
\]
Then $X^{\alpha}(F)$\ is an $\mathcal{F}$-closed $\alpha$-invariant subspace
of $X.$ Using the observation above about $sp_{\alpha}(x^{\ast})$ we have%
\[
X^{\alpha}(F)^{\ast}=X^{\alpha}(-F)
\]
where $-F=\left\{  -\gamma:\gamma\in F\right\}  .$\ If $F=\left\{
\gamma\right\}  ,\gamma\in\widehat{G}$ we will denote
\[
X^{\alpha}(\left\{  \gamma\right\}  )=X_{\gamma}.
\]
In particular for $\gamma=0$ we denote%
\[
X_{0}=X^{\alpha}%
\]
Then, we have the following consequence of the above discussions

\textbf{2.10. Corollary} \textit{Let (}$X,G,\alpha$\textit{) be an
}$\mathcal{F}$\textit{-dynamical system and }$\gamma\in\widehat{G}.$\textit{
We have\newline i) If }$x\in G$\textit{\ is such that }$sp(x)=\left\{
\gamma\right\}  ,$\textit{\ then }$\alpha_{t}(x)=\left\langle t,\gamma
\right\rangle x$\textit{ for all }$t\in G.$\textit{\newline If in addition
}$G$\textit{\ is compact we also have\newline ii) }$X_{\gamma}=\left\{
\int\overline{\left\langle t,\gamma\right\rangle }\alpha_{t}(x)dt:x\in
X\right\}  .$\textit{\newline iii) The mapping }$P_{\gamma}^{\alpha}:$\textit{
}$x\rightarrow x_{\gamma}=\int\overline{\left\langle t,\gamma\right\rangle
}\alpha_{t}(x)dt$\textit{\ from }$X$\textit{\ to }$X_{\gamma}$\textit{\ is an
}$\mathcal{F}$\textit{-continuous linear projection of }$X$\textit{\ onto
}$X_{\gamma}.$

\bigskip

\begin{proof}
Part i) follows from [17, Cor. 8.1.8.]. Part ii) is an easy consequence of the
definition of $sp(x).$\newline iii) Since obviously $x_{\gamma}=\alpha_{f}%
(x)$\ where $f(t)=\overline{\left\langle t,\gamma\right\rangle },$\ this part
follows from the discussion above.
\end{proof}

\bigskip

If $U$\ is an open subset of $\widehat{G}$, denote%
\[
X^{\alpha}(U)=\overline{lin\left\{  \alpha_{f}(x):x\in X,f\in L^{1}(G)\text{
and \textit{supp}}(\widehat{f})\subset U,\text{ compact}\right\}  }^{\sigma}.
\]
where the notation \textit{supp}$(\widehat{f})$\ stands for the support of the
Fourier transform of $f.$\ Certainly, the above spectral subspaces can be
defined similarly for the pair $(\mathcal{F}$,$X)$ with the dual
representation $\alpha^{\ast}$ of $G.$ The next lemma is a collection of known
facts about spectral subspaces. These facts are true for every locally compact
abelian group $G.$

\bigskip

In the next lemma, the notations $U,U_{\iota}$ \ stand for open subsets\ and
$F,F_{\iota}$\ for closed subsets of $\widehat{G}.$

\bigskip

\textbf{2.11. Lemma }\textit{Let (}$X,\mathcal{F}$\textit{) be a a dual pair
of Banach spaces and }$\alpha$\textit{\ a representation of }$G$\textit{\ on
}$X$\textit{. Then\newline i) }$X^{\alpha}(U)=\overline{\sum X^{\alpha
}(F_{\iota})}^{\sigma}$\textit{ if }$U=\cup F_{\iota}=\cup F_{\iota}^{\circ}%
.$\textit{\newline ii) }$X^{\alpha}(F)=\cap X^{\alpha}(U_{\iota})$\textit{ if
}$F=\cap U_{\iota}=\cap\overline{U_{\iota}}.$\textit{\newline iii) }%
$\overline{\sum X^{\alpha}(U_{\iota})}^{\sigma}=X^{\alpha}(\cup U_{\iota}%
).$\textit{\newline iv) }$X^{\alpha}(\cap F_{\iota})=\cap X^{\alpha}(F_{\iota
}).$\textit{\newline If }$X$\textit{\ is a C*-algebra and }$(X,G$%
\textit{,}$\alpha)$\textit{ an }$\mathcal{F}$\textit{-dynamical system, then
we also have\newline v) }$X^{\alpha}(F_{1})X^{\alpha}(F_{2})\subset X^{\alpha
}(\overline{F_{1}+F_{2}}),$\textit{ in particular }$X^{\alpha}(F_{1}%
)X^{\alpha}(F_{2})\subset X^{\alpha}(F_{1}+F_{2})$\textit{ if }$F_{1},F_{2}%
$\textit{\ are compact.\newline vi) }$X^{\alpha}(U_{1})X^{\alpha}%
(U_{2})\subset X^{\alpha}(U_{1}+U_{2})$\textit{ if }$U_{1},U_{2}$\textit{ are
open sets.\newline vii) Let }$x\in X$\textit{\ be such that }$sp(x)=F$%
\textit{\ is compact. Then the map }$t\rightarrow\alpha_{t}(x)$\textit{\ is
norm continuous.\newline viii) Let }$F$\textit{\ be a closed subset of
}$\widehat{G}.$\textit{ Then}%
\[
X^{\alpha}(F)=\left\{  x\in X:\alpha_{f}(x)=0\text{ if \textit{supp }}%
\widehat{f}\subset\widehat{G}-F\right\}
\]

\bigskip

\begin{proof}
Parts i), ii), iii) and iv) follow from [17, Theorem 8.1.4. viii), vii), iii)
and ii)] since, under our conditions 1)-4), the action $\alpha$\ is integrable
as noticed in [29, page 214]. Part v) follows from [29, Corollary 2.3. ii)]
for closed $F_{1},F_{2}$\ and, if $F_{1},F_{2}$\ are compact, then clearly,
$F_{1}+F_{2}$ is compact in $G,$\ so closed.\newline vi) To prove that part
vi) follows from i) and v), notice that for every open subset $U\subset
G$\ there exists a family $\left\{  F_{\iota}\right\}  $\ of compact subsets
of $G$\ such that $U=\cup F_{\iota}=\cup F_{\iota}^{\circ}.$ Hence $U_{1}=\cup
F_{1\iota}=\cup F_{1\iota}^{\circ}$\ and $U_{2}=\cup F_{2\nu}=\cup F_{2\nu
}^{\circ}$\ $\ $where $F_{1\iota}$\ and $F_{2\nu}$\ are compact. Moreover,
since $F_{1\iota}^{\circ}+F_{2\nu}^{\circ}\subset(F_{1\iota}+F_{2\nu})^{\circ
}\subset F_{1\iota}+F_{2\nu},$ we have%
\[
U_{1}+U_{2}=\cup(F_{1\iota}^{\circ}+F_{2\nu}^{\circ})=\cup(F_{1\iota}+F_{2\nu
})=\cup(F_{1\iota}+F_{2\nu})^{\circ}.
\]
and we can apply i) and v).\newline vii) By [20, Theorem 2.6.2.], there exists
a function $g\in L^{1}(G)$\ such that $\widehat{g}=1$ on an open set, $U$,
containing $sp(x).$\ Then, if $f\in L^{1}(G),$\ and \textit{supp}$\widehat{f}$
is a compact subset of $U$\ it follows that $\widehat{g}\widehat{f}%
=\widehat{g\ast f}=\widehat{f}.$ Therefore, $\alpha_{g}(\alpha_{f}%
(z))=\alpha_{f}(z)$ for all $z\in X.$ Since the set
\[
\left\{  \alpha_{f}(z):z\in X,f\in L^{1}(G),\ \text{with \textit{supp }%
}\widehat{f\text{ }}\text{ a compact subset of }U\right\}
\]
is, by definition a total subset of $X^{\alpha}(U)$\ and $F\subset U$\ it
follows, in particular, that $\alpha_{g}(x)=x.$ Therefore%
\[
\left\Vert \alpha_{t}(x)-x\right\Vert =\left\Vert \alpha_{t}(\alpha
_{g}(x))-\alpha_{g}(x)\right\Vert =\left\Vert \alpha_{g_{t}}(x)-\alpha
_{g}(x)\right\Vert \leq\left\Vert g_{t}-g\right\Vert _{1}\left\Vert
x\right\Vert .
\]
where $g_{t}$\ is the $t$-translate of $g.$ Since $\lim_{t\rightarrow
0}\left\Vert g_{t}-g\right\Vert _{1}=0,$ we are done.\newline viii) This
follows from [2, the discussion after Definition 3.2.].\newline
\end{proof}

\bigskip

Suppose that $X$\ is a C*-algebra and $(X,G,\alpha)$ an $\mathcal{F}%
$-dynamical system. We will define the Arveson and Connes spectra of the
action $\alpha.$ The Arveson spectrum of $\alpha$\ is the following subset of
the dual $\Gamma=\widehat{G}$ of $G$
\[
sp(\alpha)=\left\{  \gamma\in\Gamma:X^{\alpha}(F)\neq(0)\text{ for every
closed neighborhood }F\text{\ of }\gamma\right\}  .
\]
Clearly, $sp(\alpha)$\ is a closed subset of $\Gamma$\ and, since, as noticed
above, $X^{\alpha}(F)^{\ast}=X^{\alpha}(-F),$\ it follows that if $\gamma\in
sp(\alpha),$\ then $-\gamma\in sp(\alpha).$\ Define the strong $\mathcal{F}%
$-Arveson spectrum of the action $\alpha$
\[
\widetilde{sp}_{\mathcal{F}}(\alpha)=\left\{  \gamma\in\Gamma:\overline
{X^{\alpha}(F)^{\ast}XX^{\alpha}(F)}^{\sigma}=X\text{ for every closed
neighborhood }F\text{\ of }\gamma\right\}  .
\]
It follows immediately that $\widetilde{sp}_{\mathcal{F}}(\alpha)$\ is a
semigroup and a closed subset of $\Gamma.$

Now let $(X,G,\alpha)$ be an $\mathcal{F}$-dynamical system. We define the
strong $\mathcal{F}$-Connes spectrum of $\alpha$\ by%
\[
\widetilde{\Gamma}_{\mathcal{F}}(\alpha)=\cap\left\{  \widetilde
{sp}_{\mathcal{F}}(\alpha|_{Y}:Y\in\mathcal{H}_{\sigma}^{\alpha}(X)\right\}
.
\]
where $\mathcal{H}_{\sigma}^{\alpha}(X)$ is the set of all globally $\alpha
-$invariant $\mathcal{F}$-closed hereditary C*-subalgebras of $X.$ Using the
above observation, it follows that $\widetilde{\Gamma}_{\mathcal{F}}(\alpha
)$\ is a closed semigroup.

The next result shows that the concept of strong $\mathcal{F}$-Connes spectrum
defined above coincides with the Connes spectrum for W*-dynamical systems (i.e
if $\mathcal{F}=X_{\ast}),$ [3] and with the strong Connes spectrum for
C*-dynamical systems (i.e if $\mathcal{F}=X^{\ast})$ as defined by Kishimoto
[11]. A different concept of Connes spectrum that is not related to the
current paper has been defined for C*-dynamical systems by Olesen [15].

\bigskip

\textbf{2.12. Proposition} \textit{i) If }$X$\textit{\ is a von Neumann
algebra and }$\mathcal{F}$\textit{=}$X_{\ast}$\textit{, then }$\widetilde
{\Gamma}_{\mathcal{F}}(\alpha)$\textit{ coincides with the Connes spectrum as
defined by Connes [3].\newline ii) If }$X$\textit{ is a C*-algebra and
}$\mathcal{F}$\textit{=}$X^{\ast},$\textit{ then }$\widetilde{\Gamma
}_{\mathcal{F}}(\alpha)$\textit{ equals the strong Connes spectrum as defined
by Kishimoto [11].\newline}

\bigskip

\begin{proof}
i) Recall that the definition of the Connes spectrum for a W*-dynamical system
$(X,G,\alpha)$ [3] is the following:%
\[
\Gamma(\alpha)=\cap_{p\in X^{\alpha}}sp(\alpha|_{pXp})\text{, where }p\text{
is an }\alpha-\text{invariant projection.}%
\]
Notice that in the case of W*-dynamical systems every $\alpha$-invariant
$\mathcal{F}$-closed hereditary subalgebra of $X$\ is of the form $pXp$ with
$p\in X^{\alpha},$ projection. Clearly, $\widetilde{\Gamma}_{\mathcal{F}%
}(\alpha)\subset\Gamma(\alpha).$ Now let $\gamma\in\Gamma=\widehat{G}.$
Suppose that $\gamma\notin\Gamma_{\mathcal{F}}(\alpha).$ Then, there exists an
$\alpha$-invariant projection $p\in X^{\alpha}$ and a closed neighborhood,
$F,$ of $\gamma$ such that $\overline{pX^{\alpha}(F)pXpX^{\alpha}(F)^{\ast}%
p}^{\sigma}\neq pXp.$ Since $\overline{pX^{\alpha}(F)pXpX^{\alpha}(F)^{\ast}%
p}^{\sigma}$ is an $\alpha$-invariant $w^{\ast}$-closed hereditary
W*-subalgebra$,$ there exists an $\alpha$-invariant projection $q\in X,$
$q<p$\ such that%
\[
\overline{pX^{\alpha}(F)pXpX^{\alpha}(F)^{\ast}p}^{\sigma}=qXq.
\]
It follows that ($p-q)X^{\alpha}(F)(p-q)X(p-q)X^{\alpha}(F)^{\ast}%
(p-q)=(0),$\ so $\gamma\notin sp(\alpha|_{(p-q)X(p-q)}).$ Hence $\gamma
\notin\Gamma(\alpha).$ Therefore, $\Gamma(\alpha)\subset\widetilde{\Gamma
}_{\mathcal{F}}(\alpha)$ and part i) is proven.\newline ii) Suppose that
$X$\ is a C*-algebra and $\mathcal{F=}X^{\ast}.$\ In this case, by Remark 2.3.
i) the set of all $\mathcal{F}$-closed hereditary C* subalgebras of $X$
coincides with the set of all norm closed hereditary C* subalgebras of $X$ and
and our claim follows from the definition of the strong Connes spectrum for
C*-dynamical systems [11, Section 2].
\end{proof}

\bigskip

\section{\bigskip$\mathcal{F}$-{\protect\LARGE dynamical systems associated
with compact abelian groups}}

Let $(X,G,\alpha)$\ be an $\mathcal{F}$-dynamical system with $G$ compact
abelian. The next results (Theorems 3.2. and 3.4.) that will be used in our
study of maximality (Section 4), extend some results of Connes [3, Prop.
2.2.2. b) and Theorem 2.4.1.] for compact abelian groups to the case of
$\mathcal{F}$-dynamical systems, as well as the result of Olesen, Pedersen and
Stormer [16, Theorem 2, i)$\Leftrightarrow$ii) ].

\bigskip

\textbf{3.1. Lemma }L\textit{et }$(X,G$\textit{,}$\alpha)$\textit{ be an
}$\mathcal{F}-$\textit{dynamical system with }$G$\textit{ compact abelian}%
$.$\textit{ If }$J$\textit{ is a two sided ideal of }$X^{\alpha}%
,$\textit{\ and }$\gamma_{0}\in\widehat{G},$\textit{\ then }%
\[
(\overline{XJX}^{\sigma})_{\gamma_{0}}=\overline{lin\left\{  X_{\gamma_{1}%
}JX_{\gamma_{2}}:\gamma_{1},\gamma_{2}\in sp(\alpha)\text{ and }\gamma
_{1}\gamma_{2}=\gamma_{0}\right\}  }^{\sigma}%
\]
\textit{In particular}%
\[
(\overline{XJX}^{\sigma})^{\alpha}=\overline{lin\left\{  X_{\gamma}JX_{\gamma
}^{\ast}:\gamma\in sp(\alpha)\right\}  }^{\sigma}%
\]

\begin{proof}
Clearly,
\[
lin\left\{  X_{\gamma_{1}}JX_{\gamma_{2}}:\gamma_{1},\gamma_{2}\in
sp(\alpha),\gamma_{1}\gamma_{2}=\gamma_{0}\right\}  \subset(XJX)_{\gamma_{0}%
}\subset(\overline{XJX}^{\sigma})_{\gamma_{0}}%
\]
and therefore,
\[
\overline{lin\left\{  X_{\gamma_{1}}JX_{\gamma_{2}}:\gamma_{1},\gamma_{2}\in
sp(\alpha),\gamma_{1}\gamma_{2}=\gamma_{0}\right\}  }^{\sigma}\subset
(\overline{XJX}^{\sigma})_{\gamma_{0}}%
\]
If $x\in X,$ by Lemma 2.11. i), there is a net $\left\{  x_{\iota}\right\}
\subset lin\left\{  X_{\gamma}|\gamma\in sp(\alpha)\right\}  $ such that
$x_{\iota}\rightarrow x$ in the $\mathcal{F}$-topology. Now, let
$x_{\gamma_{1}}\in X_{\gamma_{1}}$ for some $\gamma_{1}\in sp(\alpha).$
Clearly, since by Definition 2.1. 5) the multiplication is separately
$\mathcal{F}$-continuous, it follows that $x_{\iota}jx_{\gamma_{1}}\rightarrow
xjx_{\gamma_{1}}$ in the $\mathcal{F}$-topology for all $j\in J.$ By Corollary
2.10 iii) $(x_{\iota}jx_{\gamma_{1}})_{\gamma_{0}}\rightarrow(xjx_{\gamma_{1}%
})_{\gamma_{0}}$ in the $\mathcal{F}$-topology. But, it is obvious that
$(x_{\iota}jx_{\gamma_{1}})_{\gamma_{0}}=((x_{\iota})_{\gamma_{0}\gamma
_{1}^{-1}}jx_{\gamma_{1}}),$ so%
\[
(XJX_{\gamma_{1}})_{\gamma_{0}}\subset\overline{lin\left\{  X_{\gamma_{1}%
}JX_{\gamma_{2}}:\gamma_{1},\gamma_{2}\in sp(\alpha),\gamma_{1}\gamma
_{2}=\gamma_{0}\right\}  }^{\sigma}%
\]
for every $\gamma_{1}\in sp(\alpha).$ It follows that
\[
(XJ(lin\left\{  X_{\gamma}:\gamma\in sp(\alpha)\right\}  )_{\gamma_{0}}%
\subset\overline{lin\left\{  X_{\gamma_{1}}JX_{\gamma_{2}}:\gamma_{1}%
,\gamma_{2}\in sp(\alpha),\gamma_{1}\gamma_{2}=\gamma_{0}\right\}  }^{\sigma}%
\]
for every finite subset $F\subset\widehat{G}.$ Using again the separate
$\mathcal{F}$-continuity of the multplication, we get
\[
(XJX)_{\gamma_{0}}\subset\overline{lin\left\{  X_{\gamma_{1}}JX_{\gamma_{2}%
}:\gamma_{1},\gamma_{2}\in sp(\alpha),\gamma_{1}\gamma_{2}=\gamma_{0}\right\}
}^{\sigma}%
\]
Applying again Corollary 2.10. iii) we have
\[
(\overline{XJX}^{\sigma})_{\gamma_{0}}=P_{\gamma_{0}}^{\alpha}(\overline
{XJX}^{\sigma})\subset\overline{P_{\gamma_{0}}^{\alpha}(XJX)}^{\sigma
}=\overline{(XJX)_{\gamma_{0}}}^{\sigma}\subset
\]%
\[
\subset\overline{lin\left\{  X_{\gamma_{1}}JX_{\gamma_{2}}:\gamma_{1}%
,\gamma_{2}\in sp(\alpha),\gamma_{1}\gamma_{2}=\gamma_{0}\right\}  }^{\sigma}%
\]
and we are done.
\end{proof}

\bigskip

The next result is an extension of a result of Connes [3, Proposition 2.2.2.
b)] to the framework of $\mathcal{F}$-dynamical systems $(X,G,\alpha)$\ with
$G$\ compact abelian. In order to state the result we make the following
notation: if $J\subset X^{\alpha}$ is an $\mathcal{F}$-closed ideal, let
$X_{J}=\overline{JXJ}^{\sigma}.$ Then $X_{J}$ is an $\mathcal{F}$-closed
hereditary C*-subalgebra of $X$ and obviously $X_{J}$ is $\alpha$-invariant so
$X_{J}\in\mathcal{H}_{\mathcal{\sigma}}^{\alpha}(X).$

\bigskip

\textbf{3.2. Theorem} $\widetilde{\Gamma}_{\mathcal{F}}(\alpha)=\cap\left\{
\widetilde{sp}_{\mathcal{F}}(\mathcal{\alpha}|_{X_{J}}):J\subset X^{\alpha
},\mathcal{F}-\text{\textit{closed ideal}}\right\}  $\textit{.}

\bigskip

\begin{proof}
Clearly, $\widetilde{\Gamma}_{\mathcal{F}}(\alpha)\subset\cap\left\{
\widetilde{sp}_{\mathcal{F}}(\mathcal{\alpha}|_{X_{J}}),J\subset X^{\alpha
}\text{ an }\mathcal{F}-\text{closed ideal}\right\}  .$ Now let $\gamma\in
\cap\left\{  \widetilde{sp}_{\mathcal{F}}(X_{J}):J\subset X^{\alpha
},\mathcal{F}-\text{closed ideal}\right\}  $ and $Y\in\mathcal{H}%
_{\mathcal{\sigma}}^{\alpha}(X).$ We will prove that $\gamma\in\widetilde
{sp}_{\mathcal{F}}(\mathcal{\alpha}|_{Y})$\ that is $\overline{Y_{\gamma
}Y_{\gamma}^{\ast}}^{\sigma}=Y^{\alpha},$ so $\gamma\in\widetilde{\Gamma
}_{\mathcal{F}}(\alpha).$ Denote by $J$\ the following ideal of $X^{\alpha}$%
\[
J=\overline{X^{\alpha}Y^{\alpha}X^{\alpha}}^{\sigma}%
\]
\ If $j\in J^{+},$\ then since $J\subset X^{\alpha}$ and obviously
$j=j^{\frac{1}{3}}j^{\frac{1}{3}}j^{\frac{1}{3}}$\ we have $j\in JX^{\alpha
}J,$ so $J^{+}\subset JX^{\alpha}J.$ Thus, $J=JX^{\alpha}J.$\ Therefore, since
$Y\in\mathcal{H}_{\mathcal{\sigma}}^{\alpha}(X),$ so $Y^{\alpha}\in
\mathcal{H}_{\sigma}(X^{\alpha}),$\ and the multiplication is separately
$\sigma$-continuous, we have%
\[
J=JX^{\alpha}J=\overline{JX^{\alpha}J}^{\sigma}=\overline{X^{\alpha}%
(Y^{\alpha}X^{\alpha}Y^{\alpha})X^{\alpha}}^{\sigma}=\overline{X^{\alpha
}Y^{\alpha}X^{\alpha}}^{\sigma}%
\]
and%
\[
\overline{Y^{\alpha}JY^{\alpha}}^{\sigma}=\overline{Y^{\alpha}X^{\alpha
}Y^{\alpha}X^{\alpha}Y^{\alpha}}^{\sigma}=Y^{\alpha}%
\]
Using the notation above, let $X_{J}=\overline{JXJ}^{\sigma}.$ Since
$\gamma\in\widetilde{sp}_{\mathcal{F}}(\alpha|_{X_{J}}),$ it follows that%
\[
(X_{J})_{\gamma}(X_{J})_{\gamma}^{\ast}=(\overline{JX_{\gamma}J}^{\sigma
})(\overline{JX_{\gamma}^{\ast}J}^{\sigma})
\]
is $\mathcal{F}$-dense in
\[
(\overline{JXJ}^{\sigma})^{\alpha}=\overline{JX^{\alpha}J}^{\sigma}=J
\]
Therefore%
\[
\overline{(X_{J})_{\gamma}(X_{J})_{\gamma}^{\ast}}^{\sigma}=J
\]
I will prove next that $Y_{\gamma}=Y^{\alpha}Y_{\gamma}Y^{\alpha}=Y^{\alpha
}X_{\gamma}Y^{\alpha}.$ Actually, for the the proof of the theorem, we need
only the obvious equality $Y_{\gamma}=\overline{Y^{\alpha}X_{\gamma}Y^{\alpha
}}^{\sigma},$ but I believe that the equality without the closure is worth
proving$.$ Clearly
\[
Y^{\alpha}X_{\gamma}Y^{\alpha}\subset Y_{\gamma}%
\]
Now let $x\in Y_{\gamma}.$ Denote $a=x^{\ast}x\in Y^{\alpha},b=xx^{\ast}\in
Y^{\alpha}.$\ Applying [17, Proposition 1.4.5.] and its proof, for any
$s>0,s<\frac{1}{2},$ fixed, we have%
\[
x^{\ast}=vb^{s}\text{ where }v=(\text{\textit{norm)}}\lim_{n\rightarrow\infty
}x^{\ast}(\frac{1}{n}+b)^{-\frac{1}{2}}b^{1-s}%
\]
so%
\[
x=b^{s}v^{\ast}=b^{s}(\text{\textit{norm)}}\lim_{n\rightarrow\infty}%
b^{1-s}(\frac{1}{n}+b)^{-\frac{1}{2}}x
\]
and%
\[
x=ua^{s}\text{ where }u=(\text{\textit{norm)}}\lim_{n\rightarrow\infty}%
x(\frac{1}{n}+a)^{-\frac{1}{2}}a^{1-s}%
\]
Therefore,%
\[
x=b^{s}[(\text{\textit{norm)}}\lim_{n\rightarrow\infty}(b^{1-s}(\frac{1}%
{n}+b)^{-\frac{1}{2}}x(\frac{1}{n}+a)^{-\frac{1}{2}}a^{1-s})]a^{s}%
\]
and, since%
\[
(\text{\textit{norm)}}\lim_{n\rightarrow\infty}(b^{1-s}(\frac{1}{n}%
+b)^{-\frac{1}{2}}x(\frac{1}{n}+a)^{-\frac{1}{2}}a^{1-s})\in Y_{\gamma}%
\]
the equality is proven.\ So, in particular,it follows that%
\[
\overline{Y^{\alpha}Y_{\gamma}Y^{\alpha}}^{\sigma}=\overline{Y^{\alpha
}X_{\gamma}Y^{\alpha}}^{\sigma}=Y_{\gamma}%
\]
Hence%
\[
\overline{JX_{\gamma}J}^{\sigma}=\overline{JY_{\gamma}J}^{\sigma}%
\]
Therefore%
\[
\overline{JX_{\gamma}J}^{\sigma}\overline{JX_{\gamma}^{\ast}J}^{\sigma
}=\overline{JY_{\gamma}J}^{\sigma}\overline{JY_{\gamma}^{\ast}J}^{\sigma
}\subset\overline{X^{\alpha}(Y_{\gamma}Y_{\gamma}^{\ast})X^{\alpha}}^{\sigma}%
\]
and thus%
\[
\overline{(X_{J})_{\gamma}(X_{J})_{\gamma}^{\ast}}^{\sigma}=\overline
{X^{\alpha}(Y_{\gamma}Y_{\gamma}^{\ast})X^{\alpha}}^{\sigma}%
\]
By multiplying the above equality by $Y^{\alpha}$\ to the left and to the
right, and taking into account the separate $\mathcal{F}$-continuity of the
multiplication, we get%
\[
\overline{Y_{\gamma}Y_{\gamma}^{\ast}}^{\sigma}=Y^{\alpha}%
\]
Therefore, $\gamma\in\widetilde{sp}_{\mathcal{F}}(\alpha|_{Y})$ and we are done.
\end{proof}

\bigskip

Let $(X,\mathcal{F})$ be a dual pair of Banach spaces with $X$ a C*-algebra.
$X$ is said to be $\mathcal{F}$-simple if $X$ does not contain any non-trivial
$\mathcal{F}$-closed two sided ideal.

\bigskip

\textbf{3.3. Remark }\textit{If }$X$\textit{\ is }$\mathcal{F}$\textit{-simple
and }$Y\subset X$\textit{\ is an }$\mathcal{F}$\textit{-closed hereditary
C*-subalgebra, then }$Y$\textit{\ is }$\mathcal{F}$\textit{-simple.}

\bigskip

\begin{proof}
Indeed, let $Y\subset X$\ be an $\mathcal{F}$-closed hereditary C*-subalgebra.
Then, it follows that
\[
Y=\overline{YXY}^{\sigma}%
\]
Indeed, since, in particular, $Y$\ is a hereditary C*-subalgebra of $X$\ it
follows that $Y=\overline{YXY.}^{\left\Vert {}\right\Vert }$\ Hence%
\[
YXY\subset Y=\overline{YXY}^{\left\Vert {}\right\Vert }\subset\overline
{YXY.}^{\sigma}%
\]
Since $Y$\ is $\mathcal{F}$-closed the claim follows. Now, if $I\subset Y$ is
an $\mathcal{F}$-closed ideal of $Y$\ then $XIX$ is an ideal of $X$. Since $X$
is $\mathcal{F}$-simple, we have $\overline{XIX}^{\sigma}=X.$ Therefore
\[
\overline{YXIXY}^{\sigma}=Y
\]
On the other hand, since $I$\ is an ideal of $Y$%
\[
\overline{YXIXY}^{\sigma}=\overline{YXYIYXY}^{\sigma}=I
\]
Thus $I=Y$ and $Y$ is $\mathcal{F}$-simple.\ 
\end{proof}

\bigskip

Let $(X,G$,$\alpha)$ be an $\mathcal{F-}$dynamical system. $X$ is called
$\alpha$-simple if $X$ does not contain any non-trivial $\mathcal{F}$-closed
$\alpha$-invariant two sided ideal. The following result extends [3,
Th\'{e}or\`{e}me 3.4.1.] and [16,Theorem 2, i)$\Leftrightarrow$ii)] to the
case of $\mathcal{F}$-dynamical systems. It will be used in the
characterization of maximality for periodic $\mathcal{F}$-dynamical systems
(Proposition 4.5. below).

\bigskip

\textbf{3.4. Theorem }L\textit{et }$(X,G$\textit{,}$\alpha)$\textit{ be an
}$\mathcal{F}-$\textit{dynamical system with }$G$\textit{ compact abelian.
Then the following are equivalent:\newline i) }$X^{\alpha}$\textit{ is
}$\mathcal{F}-$\textit{simple\newline ii) }$X$\textit{ is }$\alpha
$\textit{-simple and }$\widetilde{\Gamma}_{\mathcal{F}}(\alpha)=sp(\alpha).$

\bigskip

\begin{proof}
i)$\Longrightarrow$ii) Suppose that $X^{\alpha}$\ is $\mathcal{F}$-simple.
Then it is immediate from the definitions that $\widetilde{sp}_{\mathcal{F}%
}(\alpha)=sp(\alpha).$ Let $\gamma\in sp(\alpha)$ be arbitrary. Applying
Theorem 3.2. it follows that $\gamma\in\widetilde{\Gamma}_{\mathcal{F}}%
(\alpha),$ so $\widetilde{\Gamma}_{\mathcal{F}}(\alpha)=sp(\alpha).$ Let us
prove that $X$ is $\mathcal{\alpha}$-simple. If $J$\ is an $\mathcal{F}%
$-closed $\alpha$-invariant ideal of $X,$\ then $J^{\alpha}$ is an ideal of
$X^{\alpha}.$ Since $X^{\alpha}$\ is $\mathcal{F}$-simple, it follows that
$J^{\alpha}=X^{\alpha},$\ so $J=X.$\newline ii)$\Longrightarrow$i) Suppose
that $X$\ is $\mathcal{\alpha}-$simple and that $\widetilde{\Gamma
}_{\mathcal{F}}(\alpha)=sp(\alpha).$ Let $J$\ be an $\mathcal{F}$-closed ideal
of $X^{\alpha}.$ We will prove first that, for every $\gamma\in sp(\alpha
)$\ we have $X_{\gamma}JX_{\gamma}^{\ast}\subset J$\ . Then, by using Lemma
3.1. we have that ($\overline{XJX}^{\sigma})^{\alpha}\subset J$ and, since
$X$\ is assumed to be $\mathcal{\alpha}$-simple it will follow $J=X^{\alpha}.$
Let now $\gamma\in sp(\alpha)=\Gamma_{\mathcal{F}}(\alpha).$ If we denote
$Y=\overline{JXJ}^{\sigma}$\ then, since $\overline{JXJ}^{\left\Vert
{}\right\Vert }$\ is a C*-hereditary subalgebra of $X,$ from Lemma 2.8.\ it
follows that $Y\in\mathcal{H}_{\mathcal{\sigma}}^{\alpha}(X).$ It is immediate
that $Y_{\gamma}=\overline{JX_{\gamma}J}^{\sigma}.$ Since $\gamma\in
\Gamma_{\mathcal{F}}(\alpha),$ we have%
\[
Y^{\alpha}=\overline{Y_{\gamma}Y_{\gamma}^{\ast}}^{\sigma}=\overline
{JX_{\gamma}JX_{\gamma}^{\ast}J}^{\sigma}%
\]
It follows that
\[
\overline{X_{\gamma}^{\ast}Y^{\alpha}X_{\gamma}}^{\sigma}=\overline{X_{\gamma
}^{\ast}JX_{\gamma}JX_{\gamma}^{\ast}JX_{\gamma}}^{\sigma}\subset
\overline{X^{\alpha}JX^{\alpha}}^{\sigma}=J
\]
Therefore
\[
X_{\gamma}^{\ast}JX_{\gamma}\subset J
\]
for every $\gamma\in sp(\alpha).$\ By Lemma 3.1.
\[
(\overline{XJX}^{\sigma})^{\alpha}=\overline{lin\left\{  X_{\gamma}JX_{\gamma
}^{\ast}|\gamma\in sp(\alpha)\right\}  }^{\sigma}%
\]
Hence%
\[
(\overline{XJX}^{\sigma})^{\alpha}\subset J
\]
Since $X$\ is $\mathcal{\alpha}-$simple, it follows that $J=X^{\alpha}$, so
$X^{\alpha}$\ is $\mathcal{F}-$simple.
\end{proof}

\bigskip

\section{Maximality of subalgebras of analytic elements associated with
one-parameter $\mathcal{F}$-dynamical systems}

This section contains our main result. Let $\left(  X,%
\mathbb{R}
,\alpha\right)  $ be a non trivial one-parameter $\mathcal{F}-$dynamical
system. Then, according to Lemma 2.11. v), $X^{\alpha}([0,\infty))$\ is an
$\mathcal{F}$-closed subalgebra of $X.$ We will prove that the maximality of
$X^{\alpha}([0,\infty))$\ among all $\mathcal{F}$-closed subalgebras of
$X$\ is characterized by a spectral condition (Condition 4.3. below). Our
result contains all the special cases considered previously ([27], [10], [7],
[13], [14], [18], [19], [24], [25], [12]). In addition our result contains
also the case of multiplier dynamical systems that was not considered before.
In particular, as we mentioned in the Introduction, if $C_{b}(%
\mathbb{R}
)$ is the C*-algebra of bounded continuous functions\ on $%
\mathbb{R}
$ and $\alpha_{t}(f)(s)=f(s-t)$ is the translation of $f$\ by $t\in%
\mathbb{R}
,$ then $H^{\infty}(%
\mathbb{R}
)\cap C_{b}(%
\mathbb{R}
)$ is a maximal strictly closed subalgebra of $C_{b}(%
\mathbb{R}
).$ I mention that the methods used by Solel [24], and [25] are specific to
von Neumann algebras and those used \ by Kishimoto [12] are mostly C*-algebra
techniques, including the use of irreducible representations, specific to C*-algebras.

Let \textit{(}$X,\mathcal{F})$\textit{ }be a dual pair of Banach spaces such
that $X$ is a C*-algebra.\textit{ }Throughout this section we will assume that
Definition 2.1., 1)-5) is satisfied. Examples of such dual pairs were given in
Remark 2.2.

\bigskip

\textbf{4.1. Proposition. }\textit{Suppose that (}$X,\mathcal{F})$\textit{ is
a dual pair of Banach spaces such that }$X$\textit{ is a C*-algebra, so
Definition 2.1., 1)-5) holds. Then the dual }$\mathcal{F}^{\ast}$\textit{ of
}$\mathcal{F}$\textit{\ is a von Neumann algebra. Actually, there exists a
central projection }$q\in X^{\ast\ast}$\textit{\ such that }$\mathcal{F}%
^{\ast}=qX^{\ast\ast}.$

\bigskip

\begin{proof}
We will prove first that $\mathcal{F}$\ is closed to involution and to
translations with elements of $X.$ Indeed let $\varphi\in\mathcal{F}.$\ If
$\left(  x_{\alpha}\right)  \subset X,$ is such that $x_{\alpha}\rightarrow
x_{0}\in X$ in the $\mathcal{F}$-topology, then since the involution in
$X$\ is $\mathcal{F}$-continuous, it follows that $x_{\alpha}^{\ast
}\rightarrow x_{0}^{\ast}$\ in the $\mathcal{F}$-topology of $X,$ so
$\varphi(x_{\alpha}^{\ast})\rightarrow\varphi(x_{0}^{\ast})$. Hence
$\varphi^{\ast}(x_{\alpha})=\overline{\varphi(x_{\alpha}^{\ast})}%
\rightarrow\overline{\varphi(x_{0}^{\ast})}=\varphi^{\ast}(x_{0})$ and
therefore $\varphi^{\ast}(x_{\alpha})\rightarrow\varphi^{\ast}(x_{0}).$\ It
follows that $\varphi^{\ast}$ is an $\mathcal{F}-$continuous\ functional on
$X.$\ Similarly, if $y\in X,$\ since the multiplication in $X$\ is separately
$\mathcal{F}$ -continuous, it follows that $(L_{y}\varphi)(x_{\alpha}%
)=\varphi(yx_{\alpha})\rightarrow\varphi(yx_{0})=(L_{y}\varphi)(x_{0}),$ so
$L_{y}\varphi$\ is an $\mathcal{F}$-continuous linear functional on $X.$Then,
according to a well known fact about dual pairs of topological vector spaces
([22, Chapter IV, 1.2.]), we have that $\varphi^{\ast}\in\mathcal{F}$\ and
$L_{y}\varphi\in\mathcal{F}.$\ Using the facts proven above, one can
immediately infer that\ the annihilator $\mathcal{F}^{\perp}$ of $\mathcal{F}%
$\ in $X^{\ast\ast}$ is a two sided ultraweakly closed ideal of $X^{\ast\ast
},$ so $\mathcal{F}^{\perp}=pX^{\ast\ast}$ for some central projection $p\in
X^{\ast\ast}$\ and therefore the dual $\mathcal{F}^{\ast}$\ of $\mathcal{F}%
,$\ which is isomorphic with the quotient $X^{\ast\ast}/\mathcal{F}^{\perp
}=(1-p)X^{\ast\ast},$ is a von Neumann algebra with identity $q=1-p$.
\end{proof}

\bigskip

We can state the following

\bigskip

\textbf{4.2. Remark. }Suppose\textbf{ }\textit{that (}$X,\mathcal{F})$\textit{
is a dual pair of Banach spaces such that }$X$\textit{\ is a C*-algebra. Then}

\textit{i) }$\overline{X}^{w}=\mathcal{F}^{\ast}$

\textit{ii) If }$B\in\mathcal{H}_{\sigma}^{\alpha}(X),$ then there exists a
projection $p\in\mathcal{F}^{\ast}$ such that $\overline{B}^{w}=(p\mathcal{F}%
^{\ast}p),$\ so $B=(p\mathcal{F}^{\ast}p)\cap X$ and, in particular,\newline%
\ \ \ \ \bigskip\textit{iii) If }$Y$\textit{\ is an }$\mathcal{F}%
$\textit{-closed ideal of }$X,$\textit{\ then there exists a central
projection }$q\in\mathcal{F}^{\ast}$\textit{\ such that }$Y=(q\mathcal{F}%
^{\ast}q)\cap X,$ so $\overline{Y}^{w}=q\mathcal{F}^{\ast}q$.\newline%
\ \ \ \textit{\ iv) If }$p\in\mathcal{F}^{\ast}$\textit{\ is a projection
which is the strong limit of an increasing net }$\left\{  e_{\lambda}\right\}
$\textit{\ of positive elements of }$X,$\textit{\ then }$pF^{\ast}p\cap
X$\textit{ is an }$\mathcal{F}$\textit{-closed hereditary C*-subalgebra of
}$X$\textit{ that contains }$\left\{  e_{\lambda}\right\}  .$\newline%
\ \ \ \textit{v) Let }$\left(  X,%
\mathbb{R}
,\alpha\right)  $\textit{ be a one- parameter }$\mathcal{F}$\textit{-
dynamical system and }$\delta>0.$\textit{\ Let }$\left\{  e_{\lambda}\right\}
_{\lambda}$\textit{\ be an approximate identity of }$X$\textit{ and }$f\in
L^{1}(%
\mathbb{R}
),f\geqslant0$\textit{\ such that\ }%
\[
\int_{G}f(p)dp=1\text{ \textit{and supp}}\widehat{f}\subset(-\delta,\delta).
\]
\textit{Then }$\left\{  d_{\lambda}\right\}  _{\lambda}$\textit{\ where }%
\[
d_{\lambda}=\alpha_{f}(e_{\lambda})=\int_{G}f(p)\alpha_{p}(e_{\lambda})dp\in
X^{\alpha}((-\delta,\delta)).
\]
\textit{is an approximate identity\ of }$X.$

\begin{proof}
Part i) follows from well known facts about dual pairs of topological vector
spaces [22, Chapter IV, 1.3.]$.$ To prove ii) notice that from Lemma 2.8
applied to the pair $(\mathcal{F}^{\ast},\mathcal{F})$\ it follows that
$\overline{B}^{w}$\ is an $\mathcal{F}$-closed hereditary subalgebra of the
von Neumann algebra $\mathcal{F}^{\ast}$ and therefore there exists a
projection $p\in\mathcal{F}^{\ast}\ $such that $\overline{B}^{w}%
=p\mathcal{F}^{\ast}p.$\ Since the pair $(X,\mathcal{F})$ is a dual pair of
Banach spaces, the Hahn Banach separation theorem implies $\overline{B}%
^{w}\cap X=B.$ Part iii) is a consequence of ii). The proof of iv) is
straightforward. Next, we will prove v). Such a function, $f,$\ exists (for
example a Fejer kernel).\ Since $\int_{%
\mathbb{R}
}f(p)dp=1$ it follows immediately that $\left\Vert d_{\lambda}\right\Vert
=\sup\left\{  \left\vert \varphi(d_{\lambda})\right\vert :\varphi
\in\mathcal{F},\left\Vert \varphi\right\Vert \leq1\right\}  \leq\int
f(p)\left\Vert \alpha_{p}(e_{\lambda})\right\Vert \leq1,$ for all $\lambda.$
Let $x\in X$\ and $p\in G.$ Since $\left\{  e_{\lambda}\right\}  _{\lambda}%
$\ is an approximate identity of $X,$ we have, for every $p\in G$%
\[
(norm)\lim_{\lambda}x\alpha_{p}(e_{\lambda})=(norm)\lim_{\lambda}\alpha
_{p}(\alpha_{p^{-1}}(x)e_{\lambda})=x.
\]
and%
\[
(norm)\lim_{\lambda}\alpha_{p}(e_{\lambda})x=(norm)\lim_{\lambda}\alpha
_{p}(e_{\lambda}\alpha_{p^{-1}}(x))=x
\]
Using the definition of $\alpha_{f},$ the assumption that $f\geqslant0,$\ the
Lebesgue dominated convergence theorem and taking into account that $\int
_{G}f(p)dp=1$ the result follows.
\end{proof}

In the rest of this paper we will assume that the $\mathcal{F}$-dynamical
system $(X,%
\mathbb{R}
,\alpha)$\ is non trivial, that is $sp(\alpha)\neq\left\{  0\right\}  $.

The next Condition will be proven to be necessary and sufficient for
maximality of subalgebras of analytic elements of non trivial $\mathcal{F}%
$-dynamical systems..

\textbf{4.3. Spectral Condition} \textbf{for one-parameter }$\mathcal{F}%
$\textbf{-dynamical systems} \textit{Let }$\left(  X,\mathbf{%
\mathbb{R}
},\alpha\right)  $\textit{ be a non trivial }$\mathcal{F}$\textit{-dynamical
system. Suppose that\newline a) Either }$sp(\alpha)=\left\{  -\gamma
_{0},0,\gamma_{0}\right\}  $\textit{ for some }$\gamma_{0}>0$\textit{ and
}$J=\overline{X_{\gamma_{0}}X_{\gamma_{0}}^{\ast}}^{\sigma}$\textit{ is an
}$\mathcal{F}$\textit{-simple ideal of }$X^{\alpha}.$\textit{\newline
or\newline b) There exists an }$\mathcal{F}$\textit{-closed }$\alpha
$\textit{-simple two sided ideal }$Y\subset X$\textit{\ such that\newline\ b1)
}$Y^{\alpha}((0,\infty))=X^{\alpha}((0,\infty))$\textit{, and so }$Y^{\alpha
}((-\infty,0))=X^{\alpha}((-\infty,0))$\textit{ as well.\newline\ b2)
}$sp(\alpha|_{Y})=\widetilde{\Gamma}_{\mathcal{F}}(\alpha|_{Y}).$%
\textit{\newline\ b3) }$\overline{Y+X^{\alpha}([0,\infty))}^{\sigma}%
$\textit{=}$X.$

\bigskip

Suppose that $\left(  X,\mathbf{%
\mathbb{R}
},\alpha\right)  $ is a periodic $\mathcal{F}$-dynamical system, that is
$sp(\alpha)\subset%
\mathbb{Z}
\gamma_{1}$\ for some $\gamma_{1}>0.$ Then the quotient group $%
\mathbb{R}
/sp(\alpha)_{\perp}=G$ is a compact group and $\widehat{G}=%
\mathbb{Z}
\gamma_{1}$. We will prove that for periodic $\mathcal{F}$-dynamical systems
the Condition 4.3. is equivalent to the following condition:

\bigskip

(\textbf{S}) There exists an $\mathcal{F}$-simple ideal $J\subset X^{\alpha}%
$\ such that for every $\gamma\in sp(\alpha),$ $\gamma>0,\overline{X_{\gamma
}X_{\gamma}^{\ast}}^{\sigma}=J.$

\bigskip

Condition (\textbf{S}) above was considered by us [19] in the particular case
of periodic C*-dynamical systems. The proof of the equivalence of the
Condition 4.3. and condition (\textbf{S}) for periodic $\mathcal{F}$-dynamical
systems will be given below. This equivalence is necessary in proving that if
$(X,%
\mathbb{R}
,\alpha)$\ is a periodic $\mathcal{F}$-dynamical system such that $X^{\alpha
}([0,\infty))$\textit{ }is a maximal $\mathcal{F}$-closed subalgebra of $X,$
then Condition 4.3. is satisfied (Lemma 4.10. below).

Let $(X,%
\mathbb{R}
,\alpha)$ be a one-parameter $\mathcal{F}$-dynamical system.\ In what follows,
we will need the following notations: If there exists $\gamma_{1}\in
sp(\alpha),\gamma_{1}>0$ such that $sp(\alpha)\subset%
\mathbb{Z}
\gamma_{1},$ then the annihilator $(sp(\alpha))_{\perp}\subset%
\mathbb{R}
$ is a discrete subgroup of $%
\mathbb{R}
$ that equals $(%
\mathbb{Z}
\gamma_{1})_{\perp}=%
\mathbb{Z}
\frac{2\pi}{\gamma_{1}}.$\ In this case we can consider the $\mathcal{F}%
$-dynamical system $(X,G,\widetilde{\alpha})$ with the compact group $G=%
\mathbb{R}
/(%
\mathbb{Z}
\gamma_{1})_{\perp}$. If $h:%
\mathbb{R}
\rightarrow G$ is the canonical mapping, then $\widetilde{\alpha}%
_{h(t)}(x)=\alpha_{t}(x).$ Clearly, $sp(\widetilde{\alpha})=sp(\alpha)\subset%
\mathbb{Z}
\gamma_{1}.$ For every $\gamma\in sp(\widetilde{\alpha})=sp(\alpha)$ let
$P_{\gamma}^{\alpha}:X\rightarrow X_{\gamma}$ be the mapping defined above%
\[
P_{\gamma}^{\alpha}(x)=\int_{G}\overline{\left\langle t,\gamma\right\rangle
}\widetilde{\alpha_{t}}(x)dt
\]
where $dt$ is the normalized Haar measure on $G$ and the integral is taken in
the $\mathcal{F}$-topology as discussed above. The next lemma is an extension
of [19, Proposition 10] to the case of $\mathcal{F}$-dynamical systems$.$

\bigskip

\textbf{4.4. Lemma} \textit{Suppose that }$(X,%
\mathbb{R}
,\alpha)$\textit{\ is a periodic }$\mathcal{F}$\textit{-dynamical system and
Condition (\textbf{S}) holds. Then, either there exists }$\gamma_{1}%
>0,$\textit{ such that }$sp(\alpha)=\left\{  -\gamma_{1},0,\gamma_{1}\right\}
$\textit{ and }$\overline{X_{\gamma_{1}}X_{\gamma_{1}}^{\ast}}^{\sigma}%
$\textit{ is a simple ideal of }$X^{\alpha},$\textit{ or there exists }%
$\gamma_{1}>0$\textit{ such that }$sp(\alpha)=%
\mathbb{Z}
\gamma_{1},$\textit{ so, in this case, }$sp(\alpha)$\textit{ is a discrete
subgroup of }$%
\mathbb{R}
.$

\bigskip

\begin{proof}
We will prove first that, if Condition (\textbf{S}) holds and $\gamma
,\gamma^{^{\prime}}\in sp(\alpha)\cap(0,\infty),$ then $\gamma-\gamma
^{^{\prime}}\in sp(\alpha).$ From the definition of spectral subspaces it
follows that $X_{-\gamma^{^{\prime}}}\ X_{\gamma}\subset X_{\gamma
-\gamma^{^{\prime}}},$ and thus if we prove that $X_{-\gamma^{^{\prime}}%
}\ X_{\gamma}\neq\left\{  0\right\}  ,\ $then it will follow that\ $\gamma
-\gamma^{^{\prime}}\in sp(\alpha).$ Suppose to the contrary that%
\[
X_{-\gamma^{^{\prime}}}\ X_{\gamma}=\left\{  0\right\}
\]
By multiplying the above equality to the left by $X_{\gamma^{\prime}}$\ and to
the right by $X_{-\gamma}$\ we get%
\[
X_{\gamma^{^{\prime}}}X_{-\gamma^{^{\prime}}}\ X_{\gamma}X_{-\gamma}=\left\{
0\right\}
\]
Taking into account that the multiplication in $X$\ is -continuous it follows
that
\[
J=JJ=\left\{  0\right\}
\]
which contradicts the assumptions that $\gamma,\gamma^{^{\prime}}\in
sp(\alpha)$ so the claim is proven. Now suppose that $sp(\alpha)$ contains
more than three points. Let $\gamma_{1}$ be the smallest positive element of
$sp(\alpha).$ Denote $\gamma_{2}=\min\left\{  \gamma\in sp(\alpha
):\gamma>\gamma_{1}\right\}  .$ Since $sp(\alpha)$\ contains more than three
points and the system is periodic, $\gamma_{2}>\gamma_{1}.$ Then $\gamma
_{2}-\gamma_{1}\in sp(\alpha).$ From the definition of $\gamma_{1}$ and
$\gamma_{2}$ it immediately follows that $\gamma_{2}=2\gamma_{1}.$ Therefore%
\[
\overline{X_{2\gamma_{1}}X_{-2\gamma_{1}}}^{\sigma}=\overline{X_{\gamma_{1}%
}X_{-\gamma_{1}}}^{\sigma}=J
\]
By multiplying the above double equality to the left by $X_{-\gamma_{1}}$ and
to the right by $X_{\gamma_{1}}$ and taking into account that, as shown above,
$X_{-\gamma^{^{\prime}}}\ X_{\gamma}\neq\left\{  0\right\}  $ for every
$\gamma,\gamma^{^{\prime}}\in sp(\alpha),$ we get that $\left\{  0\right\}
\neq X_{-\gamma_{1}}\ X_{\gamma_{1}}\subset J.$ Since $J$\ is $\mathcal{F}%
$-simple, and clearly $X_{-\gamma_{1}}X_{\gamma_{1}}$ is an ideal, it follows
that
\[
\overline{X_{-\gamma_{1}}\ X_{\gamma_{1}}}^{\sigma}=J
\]
This equality and the definition of $J$ implies
\[
X_{\gamma_{1}}=\overline{X_{\gamma_{1}}J}^{\sigma}%
\]
and
\[
\overline{X_{\gamma_{1}}JX_{-\gamma_{1}}}^{\sigma}=J
\]
It follows that%
\[
\overline{X_{\gamma_{1}}X_{2\gamma_{1}}X_{-2\gamma_{1}}X_{-\gamma_{1}}%
}^{\sigma}=J
\]
Hence $X_{\gamma_{1}}X_{2\gamma_{1}}\neq\left\{  0\right\}  $ and thus
$X_{3\gamma_{1}}\neq\left\{  0\right\}  ,$ so $3\gamma_{1}\in sp(\alpha).$ By
induction it follows that $sp(\alpha)=%
\mathbb{Z}
\gamma_{1}$ and the lemma is proven.
\end{proof}

We can prove now the following

\bigskip

\textbf{4.5. Proposition} \textit{Let }$\left(  X,\mathbf{%
\mathbb{R}
},\alpha\right)  $\textit{ be a periodic }$\mathcal{F}$\textit{-dynamical
system. Then Condition 4.3. and Condition (\textbf{S}) are equivalent.}

\bigskip

\begin{proof}
Suppose Condition 4.3. holds. If $sp(\alpha)$\ contains only three points we
have nothing to prove. Suppose that $sp(\alpha)$\ contains more than three
points. By Condition 4.3. b2), there exists $\gamma_{1}\in%
\mathbb{R}
,\gamma_{1}>0$ such that $sp(\alpha)=%
\mathbb{Z}
\gamma_{1}.$\ Let $Y\subset X$ be the $\alpha$-simple ideal of $X$\ as in
Condition 4.3. b). Clearly, since $Y\ $is an ideal of $X,$ it follows that
$Y^{\alpha}$\ is an ideal of $X^{\alpha}.$\ Then, using Condition 4.3. b2) and
Theorem 3.4. it follows that $Y^{\alpha}$ is $\mathcal{F}$-simple. Therefore,
if $\gamma\in sp(\alpha),$ since $\overline{Y_{\gamma}Y_{\gamma}^{\ast}%
}^{\sigma}$ is an $\mathcal{F}$-closed ideal of $Y^{\alpha},$\ it follows that
$\overline{Y_{\gamma}Y_{\gamma}^{\ast}}^{\sigma}=Y^{\alpha},$ so Condition
(\textbf{S}) is satisfied.\newline Now suppose that Condition (\textbf{S})
holds true.\ If Condition 4.3., a) does not hold then, by the previous Lemma
4.4., $sp(\alpha)$\ is a discrete subgroup of $%
\mathbb{R}
.$ Then, as noticed in the proof of Lemma 4.4., $\overline{X_{\gamma}%
X_{\gamma}^{\ast}}^{\sigma}=J$ for every $\gamma\in sp(\alpha),\gamma\neq0.$
Let $Y=\overline{XJX}^{\sigma}$ and $\gamma_{0}\in sp(\alpha),\gamma_{0}%
\neq0.$ Applying Lemma 3.1. and taking into account Condition (\textbf{S}), we
get
\[
X_{\gamma_{0}}=\overline{JX_{\gamma_{0}}}^{\sigma}=\overline{X_{\gamma_{0}%
}JX_{-\gamma_{0}}X_{\gamma_{0}}}^{\sigma}\subset Y_{\gamma_{0}}\subset
X_{\gamma_{0}}.
\]
so $Y_{\gamma_{0}}=X_{\gamma_{0}}$ and therefore Condition 4.3. b1) is
satisfied. Applying again Lemma 3.1. and Condition (\textbf{S}), it follows
that $Y^{\alpha}=J.$ Since $Y^{\alpha}$\ is $\mathcal{F}$-simple, Theorem 3.4.
implies that $Y$\ is $\alpha$-simple and $\widetilde{\Gamma}_{\mathcal{F}%
}(\alpha|_{Y})=sp(\alpha|_{Y}).$\ By Lemma 2.11. i), Condition 4.3. b3) is
satisfied and we are done.
\end{proof}

\bigskip

In the next Proposition we will discuss the Condition 4.3. a). We will
describe the structure of the systems satisfying Condition 4.3. a) and will
prove that Condition 4.3. a) implies the maximality of $X^{\alpha}%
([0,\infty))$\textit{.}

\bigskip

\textbf{4.6. Proposition }\textit{Suppose that Condition 4.3. a) is satisfied.
Then}\newline\textit{i) }$I=\overline{X_{\gamma_{0}}^{\ast}X_{\gamma_{0}}%
}^{\sigma}$\textit{ is an }$\mathcal{F}$\textit{-simple ideal of }$X^{\alpha}%
$\textit{\ and }$IJ=\left\{  0\right\}  .$\textit{ \newline ii) }%
$Z=\overline{(I+J)X(I+J)}^{\sigma}$\textit{ is an }$\mathcal{F}$%
\textit{-simple C*-algebra and if }$p_{1}$\textit{\ and }$p_{2}$\textit{ \ are
the projections in }$\mathcal{F}^{\ast}$\textit{\ corresponding to the
hereditary C*-subalgebras }$\overline{JXJ}^{\sigma}$ and $\overline
{IXI}^{\sigma}$\textit{\ of }$X$ \textit{as in Remark 4.2. ii) then
}$X^{\alpha}([0,\infty))=\left\{  x\in X:p_{2}xp_{1}=0\right\}  .$%
\textit{\newline iii)} $X^{\alpha}([0,\infty))$\textit{ is a maximal
}$\mathcal{F}$\textit{-closed subalgebra of }$X$.

\bigskip

\begin{proof}
i\textit{) }Let $I_{1}$\ be an $F$-closed ideal of $I.$\ Then $X_{\gamma_{0}%
}I_{1}X_{-\gamma_{0}}$\ is an $\mathcal{F}$-dense ideal of $J.$\ Since $J$\ is
an $\mathcal{F}$-simple C*-subalgebra of $X^{\alpha}$\ it follows that
$\overline{X_{\gamma_{0}}I_{1}X_{-\gamma_{0}}}^{\sigma}=J$ and therefore, by
multiplying this equation to the left by $X_{-\gamma_{0}}$\ and to the right
by $X_{\gamma_{0}}$\ it follows that $I_{1}=\overline{II_{1}I}^{\sigma}=I$ and
therefore $I$\ is $\mathcal{F}$-simple.\newline ii) We will prove first that
$Z=(p_{1}+p_{2})\mathcal{F}^{\ast}(p_{1}+p_{2})\cap X$ is $\alpha$-simple.
Notice that $Z=\overline{(I+J)X(I+J)}^{\sigma}$ and $Z^{\alpha}=\overline
{I+J}^{\sigma}.$ Then, if $K\subset Z$ is a $\sigma$-closed $\alpha$-invariant
ideal, it follows that $K^{\alpha}$ is an ideal of $\overline{I+J}^{\sigma}.$
If $K^{\alpha}J=\left\{  0\right\}  $ then, if $x\in K_{\gamma_{0}}$ we have
$xx^{\ast}\in K^{\alpha}\cap J,$ so $x=0.$ Therefore, in this case,
$K=K^{\alpha}=I.$ But $I$\ cannot be an ideal of $Z$\ since $IX_{-\gamma_{0}%
}=X_{-\gamma_{0}}\nsubseteqq I.$ Similarly, it can be shown that the situation
$K^{\alpha}I=\left\{  0\right\}  $ cannot occur. Therefore $K^{\alpha}%
J\neq\left\{  0\right\}  $ and $K^{\alpha}I\neq\left\{  0\right\}  .$ Since
$I,J$ are $\mathcal{F}$-simple ideals of $X^{\alpha}$\ and $K^{\alpha}$\ is an
ideal of $\overline{I+J}^{\sigma}$\ it follows that $K^{\alpha}=\overline
{I+J}^{\sigma},$ so $K=Z$. Hence $Z$\ is $\alpha$-simple. We will prove next
that $Z$\ is $\mathcal{F}$-simple. Since $%
\mathbb{R}
$ is connected and $sp(\alpha)$\ is finite, by [17, 8.1.12]\ it follows that
$\alpha$\ is uniformly continuous. On the other hand, it is clear that the
spectrum of the action $\alpha^{\ast}$ induced by $\alpha$\ on $\mathcal{F}$
is also compact (and equals $\left\{  -\gamma_{0},0,\gamma_{0}\right\}  )$ as
is the spectrum of $\alpha^{\ast\ast}$\ the action on $\mathcal{F}^{\ast}$
induced by $\alpha^{\ast}.$ Therefore, $\alpha^{\ast\ast}$\ is uniformly
continuous. According to [17, 8.5.3.], there exists a uniformly continuous
unitary representation $t\rightarrow u_{t}$ of $%
\mathbb{R}
$\ into $\overline{Z}^{w}\subset\mathcal{F}^{\ast}$\ such that
\[
\alpha_{t}^{\ast\ast}(x)=u_{t}xu_{t}^{\ast}\text{ for all }t\in%
\mathbb{R}
,x\in\overline{Z}^{w}.
\]
It follows, in particular, that every central projection of $\overline{Z}^{w}%
$\ is $\alpha^{\ast\ast}$-invariant. This means that every $\mathcal{F}%
$-closed ideal of $Z$ is (globally) $\alpha$-invariant. Since by the previous
arguments $Z$\ is $\alpha$-simple, it follows that $Z$\ is $\mathcal{F}%
$-simple. The fact that $X^{\alpha}([0,\infty))=\left\{  x\in X:p_{2}%
xp_{1}=0\right\}  $ is straightforward.\newline iii) Let $W$\ be an
$\mathcal{F}$-closed subalgebra of $X$\ such that $X^{\alpha}([0,\infty
))\subset W.$ Then, by [19, Lemma 4 and Corollary 6], $W$\ is $\alpha
$-invariant. If $X^{\alpha}([0,\infty))\varsubsetneq W$, then $W_{-\gamma_{0}%
}\neq\left\{  0\right\}  .$ It is easy to check that $W_{-\gamma_{0}%
}W_{-\gamma_{0}}^{\ast}$\ is an ideal of $I.$ Therefore, since $I$\ is
$\mathcal{F}$-simple, $\overline{W_{-\gamma_{0}}W_{-\gamma_{0}}^{\ast}%
}^{\left\Vert {}\right\Vert }$ is $\mathcal{F}$-dense in $I.$ Since,
obviously, $W_{-\gamma_{0}}W_{-\gamma_{0}}^{\ast}$\ is a norm dense two sided
ideal of $\overline{W_{-\gamma_{0}}W_{-\gamma_{0}}^{\ast}}^{\left\Vert
{}\right\Vert },$\ we can apply [5,\ Proposition 1.7.2.]. Therefore, there
exists an approximate identity $\left(  e_{\lambda}\right)  $\ of
$\overline{W_{-\gamma_{0}}W_{-\gamma_{0}}^{\ast}}^{\left\Vert {}\right\Vert }%
$\ contained in $W_{-\gamma_{0}}W_{-\gamma_{0}}^{\ast}.$ It then follows that
$w^{\ast}-\lim e_{\lambda}=p_{2}$ in $\mathcal{F}^{\ast}.$ On the other hand,
since $\overline{IX_{-\gamma_{0}}}^{\sigma}=X_{-\gamma_{0}},$ it follows that
$(e_{\lambda})$ is a left approximate identity of $X_{-\gamma_{0}}$ in the
$\mathcal{F}$-topology of $X.$ So, if $x\in X_{-\gamma_{0}}$ we have
$e_{\lambda}x\in W_{-\gamma_{0}}$ and $e_{\lambda}x\rightarrow x$ in the
$\mathcal{F}$-topology of $X.$ Since the spectral subspace $W_{-\gamma_{0}}$
is $\mathcal{F}$-closed it follows that $x\in W.$ Hence $W=X$.
\end{proof}

\bigskip

Notice that the projections $p_{1}$\ and $p_{2}$\ may not belong to the
multiplier algebra of $X$\ as claimed in [12, Theorem 3.7. i)].

\bigskip

\textbf{4.7. Lemma }\textit{Let }$Y$\textit{ be a C*-algebra of operators on a
Hilbert space }$H$\textit{ whose weak closure, }$\overline{Y}^{w}$\textit{ is
a von Neumann algebra, i.e. contains the identity }$1$\textit{ of }%
$B(H),$\textit{ where }$B(H)$\textit{\ is the algebra of all linear bounded
operators on }$H$\textit{. Let }$T\subset Y$\textit{ be a norm closed subspace
of }$Y$\textit{ such that the norm-closed hereditary C*-subalgebra
}$C=\overline{TYT^{\ast}}^{\left\Vert {}\right\Vert }$\textit{ of }%
$Y$\textit{\ is }$w$\textit{-dense in }$Y$\textit{. Then there exists a
directed family }$\left\{  e_{\lambda}\right\}  \subset\overline
{\text{\textit{alg}}(TT^{\ast})}^{\left\Vert {}\right\Vert }\cap Y_{1}^{+}%
$\textit{ such that }$\lim_{\lambda}e_{\lambda}=1$\textit{ in the strong
operator topology of }$\overline{Y}^{w}$, where $\overline{\text{\textit{alg}%
}(TT^{\ast})}^{\left\Vert {}\right\Vert }$\ denotes the norm closed subalgebra
of $Y$ generated by $TT^{\ast}.$\textit{ }

\bigskip

\begin{proof}
Let $\Lambda$ be the collection of all finite subsets $\left\{  x_{i}\right\}
_{i=1}^{n}\subset T,n\in%
\mathbb{N}
.$ If $\lambda\in\Lambda,$\ set $a_{\lambda}=\sum_{i=1}^{n}x_{i}x_{i}^{\ast}$
and $e_{\lambda}=a_{\lambda}(\frac{1}{n}+a_{\lambda})^{-1}.$ It is clear that
$e_{\lambda}\in\overline{\text{\textit{alg}}(TT^{\ast})}^{\left\Vert
{}\right\Vert }\cap Y_{1}^{+}$ and by Lemma 2.9., $e_{\lambda}\in C.$ The next
arguments are inspired by [5, proof of Proposition 1.7.2], even if that result
is different from the statement of our current Lemma 4.7. One can easily check
that for each $i=1,...,n.$%
\[
\left\Vert (1-e_{\lambda})x_{i}x_{i}^{\ast}(1-e_{\lambda})\right\Vert
\leqslant\left\Vert \sum_{j=1}^{n}(1-e_{\lambda})x_{j}x_{j}^{\ast
}(1-e_{\lambda})\right\Vert =\frac{1}{n^{2}}\left\Vert a_{\lambda}(\frac{1}%
{n}+a_{\lambda})^{-2}\right\Vert .
\]
and, since as noticed in [5, Proof of Proposition 1.7.2.], $t(\frac{1}%
{n}+t)^{-2}\leq\frac{n}{4}$ for all $n$ $\in%
\mathbb{N}
,t\in%
\mathbb{R}
^{+},$ we have%
\[
\left\Vert (1-e_{\lambda})x_{i}x_{i}^{\ast}(1-e_{\lambda})\right\Vert
\leqslant\frac{1}{4n}.
\]
Therefore, $\lim_{\lambda}\left\Vert (1-e_{\lambda})x\right\Vert =0$ for all
$x\in T.$ Hence, $\lim_{\lambda}\left\Vert (1-e_{\lambda})xyz^{\ast
}\right\Vert =0$ for all $x,z\in T,y\in Y.$ It then follows that
$\lim_{\lambda}\left\Vert (1-e_{\lambda})z\right\Vert =0$ for all
$z\in\overline{TYT^{\ast}}^{\left\Vert {}\right\Vert }=C$. Then $\overline
{C}^{w}=\overline{Y}^{w}\subset B(H)$ and certainly, $1$\ is the unit of
$\overline{C}^{w}.$ Let $\xi\in H.$ Then, for every $\epsilon>0$\ there exists
$z\in C$ such that $\left\Vert z\xi-\xi\right\Vert <\frac{\epsilon}{3}.$Since
$\left\{  e_{\lambda}\right\}  $ is an approximate identity of $C,$\ there
exists $\lambda_{0}$ such that $\left\Vert e_{\lambda}z-z\right\Vert
<\frac{\epsilon}{3}$ for all $\lambda\succcurlyeq\lambda_{0}.$ Therefore%
\[
\left\Vert e_{\lambda}\xi-\xi\right\Vert \leq\left\Vert e_{\lambda}(\xi
-z\xi)\right\Vert +\left\Vert (e_{\lambda}z-z)\xi\right\Vert +\left\Vert
z\xi-\xi\right\Vert <\epsilon.
\]
for every $\lambda\succcurlyeq\lambda_{0}.$
\end{proof}

\bigskip The next lemma is an extension of [12, Lemma 3.1.] to the case of
dual $\mathcal{F}$-dynamical systems and general locally compact groups $G.$

\bigskip

\textbf{4.8. Lemma }\textit{Let }$\left(  Y,G,\alpha\right)  $\textit{ be an
}$\alpha$\textit{-simple }$\mathcal{F}$\textit{-dynamical system with }%
$G$\textit{\ a locally compact abelian group such that }$sp$\textit{(}%
$\alpha)=\widetilde{\Gamma}_{\mathcal{F}}(\alpha).$\textit{ If }$C\in
H_{\sigma}^{\alpha}(Y)$\textit{ and }$\gamma\in sp(\alpha),$\textit{ then
}$\overline{Y^{\alpha}(V)CY^{\alpha}(V)^{\ast}}^{\sigma}=Y$\textit{ for every
open neighborhood, }$V,\ $\textit{of }$\gamma.$

\bigskip

\begin{proof}
Since $sp$($\alpha)=\widetilde{\Gamma}_{\mathcal{F}}(\alpha),$\ it follows
that $\widetilde{\Gamma}_{\mathcal{F}}(\alpha)$ is a closed subgroup of
$\Gamma=\widehat{G}.$ Let $U$ be a fixed neighborhood of $0\in\Gamma$ so
$V=U+\gamma$ is a fixed neighborhood of $\gamma.$ Further, let $\left\{
\gamma_{\iota}\right\}  _{\iota\in I}$ be a dense subset of $\widetilde
{\Gamma}_{\mathcal{F}}(\alpha)$\ and, for each $\iota\in I,$ let $U_{\iota}%
\ $be an open neighborhood of $\gamma_{\iota}$ such that $U_{\iota}-U_{\iota
}\subset U.$ Since $\Gamma_{\mathcal{F}}(\alpha)\subset\cup U_{\iota},,$ by
Lemma 2.11. i) it follows that
\[
Y^{\alpha}(\widetilde{\Gamma}_{\mathcal{F}}(\alpha))=Y\subset\overline
{\sum_{\iota\in I}Y^{\alpha}(U_{\iota})}^{\sigma}.
\]
so%
\[
\overline{\sum_{\iota\in I}Y^{\alpha}(U_{\iota})}^{\sigma}=Y.
\]
Therefore%
\[
\overline{\sum_{\iota\in I}Y^{\alpha}(U_{\iota})}^{\sigma}C=YC.
\]
and%
\[
C\overline{\sum_{\iota\in I}Y^{\alpha}(U_{\iota})}^{\sigma}=CY.
\]
On the other hand, since $C\in\mathcal{H}_{\sigma}^{\alpha}(Y)$ and
$\widetilde{\Gamma}_{\mathcal{F}}(\alpha)$ is a group, for each $\iota\in
I$\ we have $\gamma_{\iota}-\gamma\in\widetilde{\Gamma}_{\mathcal{F}}%
(\alpha),$\ so $\overline{C^{\alpha}(-U_{\iota}+\gamma)CC^{\alpha}(-U_{\iota
}+\gamma)^{\ast}}^{\sigma}=C$. Therefore
\[
Y^{\alpha}(U_{\iota})C=Y^{\alpha}(U_{\iota})\overline{C^{\alpha}(-U_{\iota
}+\gamma)CC^{\alpha}(-U_{\iota}+\gamma)^{\ast}}^{\sigma}\subset
\]%
\[
\subset\overline{Y^{\alpha}(U_{\iota})C^{\alpha}(-U_{\iota}+\gamma)CC^{\alpha
}(-U_{\iota}+\gamma)^{\ast}}^{\sigma}\subset\overline{Y^{\alpha}(U+\gamma
)C}^{\sigma}.
\]
Similarly, for each $\kappa\in I$%
\[
CY^{\alpha}(U_{\kappa})\subset\overline{CY^{\alpha}(U+\gamma)^{\ast}}^{\sigma}%
\]
Since the multiplication in $Y$\ is separately $\mathcal{F}$-continuous, it
follows that for every $\iota,\kappa\in I$%
\[
Y^{\alpha}(U_{\iota})CY^{\alpha}(U_{\kappa})\subset\overline{Y^{\alpha
}(U+\gamma)CY^{\alpha}(U+\gamma)^{\ast}}^{\sigma}%
\]
so%
\[
YCY\subset\overline{Y^{\alpha}(U+\gamma)CY^{\alpha}(U+\gamma)^{\ast}}^{\sigma}%
\]
Since $Y$ is $\alpha$-simple, the proof is complete.
\end{proof}

\bigskip

\textbf{4.9. Theorem }\textit{Let (}$X,%
\mathbb{R}
,\alpha$\textit{) be a one-parameter }$\mathcal{F}$\textit{-dynamical system
satisfying Condition 4.3. Then }$X^{\alpha}([0,\infty))$\textit{ is a maximal
}$\mathcal{F}$\textit{-closed subalgebra of }$X$\textit{.}

\bigskip

\begin{proof}
If Condition 4.3. a) is satisfied, the conclusion follows from Proposition
4.6. iii). Suppose that Condition 4.3. b) is satisfied and let $Y$\ be the
ideal considered in that condition. By Condition 4.3. b3), $sp$($\alpha
|_{Y})=\widetilde{\Gamma}_{\mathcal{F}}(\alpha|_{Y})$ and therefore
$sp$($\alpha)$\ is a closed subgroup of $%
\mathbb{R}
.$\ It follows that either $sp(\alpha)$ is discrete, i.e. there exists
$\gamma\in%
\mathbb{R}
$ such that $sp(\alpha)=\left\{  n\gamma:n\in%
\mathbb{Z}
\right\}  $ or $sp$($\alpha)=%
\mathbb{R}
.$ Let $Z$ be an $\mathcal{F}$-closed subalgebra of $X$\ that contains
properly $X^{\alpha}([0,\infty)).$ According to [19, Lemma 4 and Corollary 6],
$Z$\ is $\alpha$-invariant. Since $X^{\alpha}([0,\infty))\subsetneqq Z,$ there
exists $\gamma\in sp(\alpha|_{Z})-[0,\infty).$ We will prove that $n\gamma\in
sp(\alpha|_{Z})$ for all $n\in%
\mathbb{N}
.$ This fact will be proved separately for each of the two possible cases:
\[
sp(\alpha)=%
\mathbb{R}
\text{ or }sp(\alpha)\text{ is a discrete subgroup of }%
\mathbb{R}
.
\]
Suppose first that $sp$($\alpha)=%
\mathbb{R}
,$ and let $r>0$ such that
\[
Z^{\alpha}((\gamma-r,\gamma+r))\subset X^{\alpha}((-\infty,0))=Y^{\alpha
}((-\infty,0))
\]
Then, there exists $x\in Z$ with $sp(x)\subset(\gamma-\frac{r}{3},\gamma
+\frac{r}{3}).$ Let $S=\left\{  \alpha_{t}(x):t\in%
\mathbb{R}
\right\}  .$ Consider the hereditary subalgebra $W\in\mathcal{H}_{\sigma
}^{\alpha}(Y),$ where $Y$\ is as in Condition 4.3. b)%
\[
W=\overline{SYS^{\ast}}^{\sigma}.
\]
Since $sp(\alpha)=\widetilde{\Gamma}_{\mathcal{F}}(\alpha)=%
\mathbb{R}
,$ it follows that $(0,\frac{r}{3})\cap\widetilde{\Gamma}_{\mathcal{F}}%
(\alpha)\neq\emptyset$, so
\[
Y=\overline{Y^{\alpha}((0,\frac{r}{3}))YY^{\alpha}((0,\frac{r}{3}))^{\ast}%
}^{\sigma}.
\]
Therefore%
\[
W=\overline{SY^{\alpha}((0,\frac{r}{3}))YY^{\alpha}((0,\frac{r}{3}))^{\ast
}S^{\ast}}^{\sigma}.
\]
The last equality above follows from the assumption of separate continuity of
the multiplication of $X$ (Definition 2.1. 5))$.$ Applying Lemma 4.8., it
follows that%
\[
Y=\overline{Y^{\alpha}((0,\frac{r}{3}))SY^{\alpha}((0,\frac{r}{3}))YY^{\alpha
}((0,\frac{r}{3}))^{\ast}S^{\ast}Y^{\alpha}((0,\frac{r}{3}))^{\ast}}^{\sigma
}.
\]
Clearly%
\[
Y^{\alpha}((0,\frac{r}{3}))SY^{\alpha}((0,\frac{r}{3}))\subset Z^{\alpha
}((\gamma-\frac{r}{3},\gamma+r))\subset Z^{\alpha}((\gamma-r,\gamma+r)).
\]
Therefore%
\[
Y=\overline{Z^{\alpha}((\gamma-r,\gamma+r))YZ^{\alpha}((\gamma-r,\gamma
+r))^{\ast}}^{\sigma}.
\]
We will prove now that $n\gamma\in sp(\alpha|_{Z})$ for all $n\in%
\mathbb{N}
.$ Multiplying the above to the left by $Z^{\alpha}((\gamma-r,\gamma+r))$\ and
to the right by $Z^{\alpha}((\gamma-r,\gamma+r))^{\ast}$, it follows that%
\[
\left\{  0\right\}  \neq Z^{\alpha}((\gamma-r,\gamma+r))Z^{\alpha}%
((\gamma-r,\gamma+r))\subset
\]%
\[
\subset Z^{\alpha}((2\gamma-2r,2\gamma+2r)).
\]
Since $r>0$ is arbitrarily small, we have that $2\gamma\in sp(\alpha|_{Z}).$
By induction, it follows that
\[
n\gamma\in sp(\alpha|_{Z})\text{ for all }n\in%
\mathbb{N}
.
\]
Now suppose that $sp$($\alpha)=%
\mathbb{Z}
\gamma_{0}$ for some $\gamma_{0}\in%
\mathbb{R}
$ and let $\gamma\in sp(\alpha|_{Z})-[0,\infty).$ By taking $V=\left\{
\gamma\right\}  $ in Lemma 4.8., we obtain similarly that $n\gamma\in
sp(\alpha|_{Z})$ for all $n\in%
\mathbb{N}
$. \newline By Remark 4.2. iii) there exists a central projection
$q\in\mathcal{F}^{\ast}$ such that $\overline{Y}^{w}=q\mathcal{F}^{\ast}$. Let
$y\in Y$ be an element such that $sp(y)$ is compact. Then, since as shown
above, $sp(\alpha|_{Z})$ is unbounded, there exists $\gamma\in sp(\alpha
|_{Z})-[0,\infty)$ and $r>0$ such that $-\gamma+sp(y)\subset\lbrack
r,\infty),$ so
\[
sp(y)\subset\lbrack\gamma+r,\infty)
\]
$.$ Notice that, since $Z^{\alpha}((\gamma-r,\gamma+r))YZ^{\alpha}%
((\gamma-r,\gamma+r))^{\ast}$ is $\mathcal{F}$-dense in $Y,$ it follows that
$\overline{Z^{\alpha}((\gamma-r,\gamma+r))YZ^{\alpha}((\gamma-r,\gamma
+r))^{\ast}}^{w}=\overline{Y}^{w}$ and therefore
\[
C=\overline{Z^{\alpha}((\gamma-r,\gamma+r))YZ^{\alpha}((\gamma-r,\gamma
+r))^{\ast}}^{\left\Vert {}\right\Vert }.
\]
is strongly dense in $\overline{Y}^{w}=q\mathcal{F}^{\ast}.$ By Lemma 4.7.
there exists a directed set
\[
\left\{  e_{\lambda}\right\}  \subset\overline{\text{\textit{alg(}}Z^{\alpha
}((\gamma-r,\gamma+r))Z^{\alpha}((\gamma-r,\gamma+r))^{\ast}))}^{\left\Vert
{}\right\Vert }\cap Y_{1}^{+}.
\]
such that $\lim_{\lambda}e_{\lambda}=q$ in the strong operator topology of
$\overline{Y}^{w}=q\mathcal{F}^{\ast}$ \ Therefore $\lim_{\lambda}e_{\lambda
}y=y$ in the strong operator topology of $q\mathcal{F}^{\ast}.$ But, since%
\begin{align*}
e_{\lambda}  &  \in\overline{\mathit{alg(}Z^{\alpha}((\gamma-r,\gamma
+r))Z^{\alpha}((\gamma-r,\gamma+r))^{\ast}))}^{\left\Vert {}\right\Vert }=\\
&  =\overline{\mathit{alg(}Z^{\alpha}((\gamma-r,\gamma+r))Z^{\alpha}%
((-\gamma-r,-\gamma+r))}^{\left\Vert {}\right\Vert }%
\end{align*}
and, by the choice of $\gamma,r$ and $y\in Z,$
\[
Z^{\alpha}((-\gamma-r,-\gamma+r))y\subset X^{\alpha}([0,\infty))\subset Z
\]
it follows that
\[
\mathit{alg(}Z^{\alpha}((\gamma-r,\gamma+r))Z^{\alpha}((-\gamma-r,-\gamma
+r))y\subset Z
\]
Since $Z$\ is, in particular, norm closed, it follows that $e_{\lambda}y\in Z$
for all $\lambda\in\Lambda.$ Since $Z$ is $\mathcal{F}-$ closed in $X$, we
have that $y\in Z$.\ Applying Lemma 2.11. i) it follows in particular that
$Y\subset Z$. Since $X^{\alpha}([0,\infty))\subset Z,$ we have $Y+X^{\alpha
}([0,\infty))\subset Z.$ By Condition 4.3. b3), $Y+X^{\alpha}([0,\infty))$ is
$\mathcal{F}$-dense in $X.$ Hence $Z=X.$
\end{proof}

\bigskip

To prove the converse of the above Theorem 4.9. we need to prove a series of
lemmas. We start with the following

\bigskip

\textbf{4.10. Lemma }\textit{Let (}$X,%
\mathbb{R}
,\alpha$\textit{) be a one-parameter }$\mathcal{F}$\textit{-dynamical system.
If }$X^{\alpha}([0,\infty))$\textit{\ is a maximal }$\mathcal{F}%
$\textit{-closed subalgebra of }$X,$\textit{\ then there exists an
}$\mathcal{F}$\textit{-closed, }$\alpha$\textit{-simple ideal, }%
$Y,$\textit{\ of }$X$\textit{\ such that\newline i) }$Y^{\alpha}%
((0,\infty))=X^{\alpha}((0,\infty))$\textit{, and so }$Y^{\alpha}%
((-\infty,0))=X^{\alpha}((-\infty,0))$\textit{ as well.\newline ii)
}$Y+X^{\alpha}([0,\infty))$\textit{ is }$\mathcal{F}$\textit{-dense in }$X.$

\bigskip

\begin{proof}
Let $J_{1}$\ be an $\alpha$-invariant ideal of $X$\ such that $J_{1}%
\nsubseteqq X^{\alpha}.$\ Since $X^{\alpha}([0,\infty))$ is a maximal
$\mathcal{F}$-closed subalgebra of $X$\ and $\overline{J_{1}+X^{\alpha
}([0,\infty))}^{\sigma}$\ is an $\mathcal{F}$-closed subalgebra properly
containing $X^{\alpha}([0,\infty)),$\ it follows that $\overline
{J_{1}+X^{\alpha}([0,\infty))}^{\sigma}=X.$ We will prove that $J_{1}^{\alpha
}((0,\infty))=X^{\alpha}((0,\infty)).$ Indeed, let $x\in X^{\alpha}%
((-\infty,0)).\ $Then, $x=\mathcal{F}-\lim(x_{\lambda}+j_{\lambda})$ where
$x_{\lambda}\in X^{\alpha}[0,\infty)$ and $j_{\lambda}\in J_{1}.$ Now let
$f\in L^{1}(%
\mathbb{R}
)$ such that $supp$($\widehat{f})\subset(-\infty,0).$ Since $\alpha_{f}$ is
$\mathcal{F}$-continuous, it follows that $\alpha_{f}(x)=\mathcal{F}%
-\lim(\alpha_{f}(x_{\lambda}+j_{\lambda})).$ Since $x_{\lambda}\in X^{\alpha
}([0,\infty)),$ by Lemma 2.11. viii), we have $\alpha_{f}(x_{\lambda})=0$ and
so, since $J_{1}$\ is $\mathcal{F}$-closed and $\alpha$-invariant,%
\[
\alpha_{f}(x)=\mathcal{F}-\lim(\alpha_{f}(j_{\lambda}))\in J_{1}.
\]
\ \ Hence $\ X^{\alpha}((-\infty,0))\subset J_{1}^{\alpha}((-\infty,0))$ so
$J_{1}^{\alpha}((-\infty,0))=X^{\alpha}((-\infty,0))$ and consequently,
$J_{1}^{\alpha}((0,\infty))=X^{\alpha}((0,\infty)).$ Let $J$\ be the
intersection of all such $J_{1}.$ Then $J$\ is an $\mathcal{F}$-closed
$\alpha$-invariant ideal of $X$\ such that $J^{\alpha}((-\infty,0))=X^{\alpha
}((-\infty,0)).$ In order to prove the existence of the ideal $Y$ in the
statement of the lemma, notice that every $\alpha$-invariant ideal $I_{1}$\ of
$J$\ must be included in the fixed point algebra $X^{\alpha}.$ Let $I$ be the
$\mathcal{F}$-closed linear span of all the ideals of $X$ that are contained
in $X^{\alpha}.$ Denote%
\[
Y=\left\{  y\in J:yI=\left\{  0\right\}  \right\}
\]
Then%
\[
Y=\left\{  y\in J:Iy=\left\{  0\right\}  \right\}
\]
Indeed, if $iy\neq0$ for some $i\in I,y\in J$ with $yI=\left\{  0\right\}
,$\ then $iyy^{\ast}i^{\ast}\neq0$\ but this contradicts the fact that
$yI=\left\{  0\right\}  .$ Then $Y$\ is an $\mathcal{F}$-closed $\alpha
$-invariant $\alpha$-simple ideal of $X.$\ The fact that $Y$ is $\mathcal{F}%
$-closed and $\alpha$-invariant is obvious from the definition of \ $Y.$ To
prove that $Y$\ is $\alpha$-simple let $K\subset Y$ be a non-zero $\alpha
$-invariant ideal of $Y.$ By the definition of $Y$\ it follows that
$K\nsubseteq X^{\alpha}.$\ Therefore $Y\subset J\subset K$ and thus $K=Y=J.$
Clearly, $Y=J$ satisfies the conditions of the Lemma.
\end{proof}

\bigskip

We will prove first the converse of Theorem 4.9. in the special case when
$0$\ is an isolated point of $sp(\alpha)$.

\bigskip

\textbf{4.11. Lemma}\bigskip\ \textit{Let (}$X,%
\mathbb{R}
,\alpha$\textit{) be a one-parameter }$\mathcal{F}$\textit{-dynamical
system}$.$\textit{Suppose that }$X^{\alpha}([0,\infty))$\textit{\ is a maximal
}$\mathcal{F}-$\textit{closed subalgebra of }$X.$\textit{\ If }$0$\textit{ is
an isolated point of }$sp(\alpha),$\textit{ then Condition 4.3. is
satisfied.\ }

\bigskip

\begin{proof}
We will prove first that there exists $\gamma_{1}\in sp(\alpha)\ $such that
$sp(\alpha)\subset%
\mathbb{Z}
\gamma_{1}.$\ Let $\gamma_{1}=\inf\left\{  \gamma\in sp(\alpha):\gamma
>0\right\}  .$ Then, since $sp(\alpha)$\ is closed and $0$\ is isolated in
$sp(\alpha),$ we have that $\gamma_{1}>0$ and $\gamma_{1}\in sp(\alpha
)$.\ Since $0$\ is isolated in $sp(\alpha)$, we have $X^{\alpha}%
((-\epsilon,\epsilon))=X^{\alpha}$\ \ for every $0<\epsilon<\gamma_{1}.$\ We
will prove next that $\gamma_{1}$\ is an isolated point of $sp(\alpha
).$\ Suppose that there exists $\gamma_{2}\in sp(\alpha)$ such that
$\gamma_{1}<\gamma_{2}<2\gamma_{1}.$ Let $0<\epsilon<\min\left\{  \frac
{\gamma_{2}-\gamma_{1}}{4},\frac{2\gamma_{1}-\gamma_{2}}{4}\right\}  $ and
denote%
\[
J_{\gamma_{1}}^{\epsilon}=\overline{X^{\alpha}((\gamma_{1}-\epsilon,\gamma
_{1}+\epsilon))X^{\alpha}((-\gamma_{1}-\epsilon,-\gamma_{1}+\epsilon
))}^{\sigma}\subset X^{\alpha}%
\]
Then, it is immediate that $J_{\gamma_{1}}^{\epsilon}$\ is a two sided ideal
of $X^{\alpha}.$\ Consider the following $\mathcal{F}$-closed\ $\alpha
$-invariant subspace of $X$%
\[
\mathcal{M}=\overline{X^{\alpha}((-\infty,0))J_{\gamma_{1}}^{\epsilon
}+X^{\alpha}([0,\infty))}^{\sigma}.
\]
and let
\[
W=\left\{  x\in X:x\mathcal{M}\subset\mathcal{M}\right\}  .
\]
Then, $W\ $is an $\mathcal{F}$-closed $\alpha$-invariant subalgebra of $X.$
Taking into account that $0$\ is isolated in $sp(\alpha)$, and that by Lemma
2.11. i) we have%
\[
X^{\alpha}([0,\infty))=
\]%
\[
\overline{X^{\alpha}+\sum_{k\in%
\mathbb{N}
}X^{\alpha}((k\gamma_{1}-\epsilon,k\gamma_{1}+\epsilon))+\sum_{m\in%
\mathbb{N}
}X^{\alpha}((m\gamma_{1}+\frac{\epsilon}{2},(m+1)\gamma_{1}-\frac{\epsilon}%
{2}))}^{\sigma}.
\]
and%
\[
X^{\alpha}((-\infty,0))=
\]%
\[
\overline{\sum_{k\in%
\mathbb{N}
}X^{\alpha}((-k\gamma_{1}-\epsilon,-k\gamma_{1}+\epsilon))+\sum_{m\in%
\mathbb{N}
}X^{\alpha}((-(m+1)\gamma_{1}+\frac{\epsilon}{2},-m\gamma_{1}-\frac{\epsilon
}{2}))}^{\sigma}.
\]
it can be checked that $X^{\alpha}([0,\infty))\subset W.$ On the other hand,
since the ideal $J_{\gamma_{1}}^{\epsilon}$\ of $X^{\alpha}$\ contains a right
approximate identity for $X^{\alpha}((-\gamma_{1}-\epsilon,-\gamma
_{1}+\epsilon)),$ it follows that $X^{\alpha}((-\gamma_{1}-\epsilon
,-\gamma_{1}+\epsilon))\subset W.$ Therefore,%
\[
X^{\alpha}([0,\infty))\subsetneqq W.
\]
Since $X^{\alpha}([0,\infty))$ is maximal, it follows that $W=X.$ Let $y\in
X^{\alpha}((\gamma_{2}-\epsilon,\gamma_{2}+\epsilon)),y\neq0.$ Then $y^{\ast
}\in W=X.$ So, in particular,\ $y^{\ast}X^{\alpha}=y^{\ast}X^{\alpha
}((-\epsilon,\epsilon))\subset W$. By Remark 4.2. v) $X^{\alpha}%
((-\epsilon,\epsilon))$ contains a Banach algebra approximate identity of $X,$
so $y^{\ast}X^{\alpha}\neq\left\{  0\right\}  .$ Since
\[
y^{\ast}X^{\alpha}\subset X^{\alpha}((-\infty,0)).
\]
we must have%
\[
y^{\ast}X^{\alpha}\subset\ \mathcal{M}^{\alpha}((-\infty,0))=\overline
{X^{\alpha}((-\infty,0))J_{\gamma_{1}}^{\epsilon}}^{\sigma}.
\]
Therefore%
\[
y^{\ast}X^{\alpha}J_{\gamma_{1}}^{\epsilon}=y^{\ast}J_{\gamma_{1}}^{\epsilon
}\neq\left\{  0\right\}  .
\]
So, in particular
\[
y^{\ast}X^{\alpha}((\gamma_{1}-\epsilon,\gamma_{1}+\epsilon))\neq0.
\]
But, according to Lemma 2.11. vi) we have%
\[
\left\{  0\right\}  \neq y^{\ast}X^{\alpha}((\gamma_{1}-\epsilon,\gamma
_{1}+\epsilon))\subset X^{\alpha}((\gamma_{1}-\gamma_{2}-2\epsilon,\gamma
_{1}-\gamma_{2}+2\epsilon).
\]
The choice of $\gamma_{2}\ $and $\epsilon$\ imply\ \ $(\gamma_{1}-\gamma
_{2}-2\epsilon,\gamma_{1}-\gamma_{2}+2\epsilon)\subset(-\gamma_{1},0).$ This
is a contradiction with the assumption that that $0$\ is an isolated point of
$sp(\alpha)$\ and $\gamma_{1}=\inf\left\{  \gamma\in sp(\alpha):\gamma
>0\right\}  =\min\left\{  \gamma\in sp(\alpha):\gamma>0\right\}  .$ Hence
$sp(\alpha)\cap(\gamma_{1},2\gamma_{1})=\emptyset.$ Similarly, one can show
that $sp(\alpha)\cap(k\gamma_{1},(k+1)\gamma_{1})=\emptyset$ for every $k\in%
\mathbb{Z}
.$ It follows that $sp(\alpha)\subset%
\mathbb{Z}
\gamma_{1}.$ We will prove next that Condition 4.3. is satisfied. Since as
shown above, $sp(\alpha)\subset%
\mathbb{Z}
\gamma_{1}$\ where $\gamma_{1}=\min\left\{  \gamma\in sp(\alpha):\gamma
>0\right\}  ,$ from Proposition 4.5. it follows that Condition 4.3. is
equivalent in this case with Condition (\textbf{S}). We will show that
Condition (\textbf{S}) is satisfied. Let $J=\overline{X_{\gamma_{1}}%
X_{\gamma_{1}}^{\ast}}^{\sigma}$ and let $I$\ be a non zero $\mathcal{F}%
$-closed ideal of $J.$ Let%
\[
\mathcal{M=}\overline{\sum_{\gamma\in sp(\alpha),\gamma>0}X_{-\gamma
}I+X^{\alpha}([0,\infty))}^{\sigma}%
\]
Then $\mathcal{M}$\ is an $\alpha$-invariant $\mathcal{F}$-closed linear
subspace\ of $X.$ Let%
\[
W=\left\{  x\in X:x\mathcal{M}\subset M\right\}
\]
Then $W$\ is an $\mathcal{F}$-closed subalgebra of $X.$\ Clearly $X^{\alpha
}([0,\infty))\subset W.$ On the other hand, since $I\neq\left\{  0\right\}  ,$
it follows that $X_{-\gamma_{1}}I\neq\left\{  0\right\}  $\ and therefore
there exists $y\in X_{-\gamma_{1}}$ and $i_{0}\in I$\ such that $z=yi_{0}%
\neq0.$ Since $\gamma_{1}=\min\left\{  \gamma\in sp(\alpha):\gamma>0\right\}
$\ it follows that $z\in W.$ Therefore, since $z\in X-X^{\alpha}([0,\infty))$
and $X^{\alpha}([0,\infty))$ is maximal, it follows that $W=X$. Hence, in
particular $X_{-\gamma_{1}}X^{\alpha}\subset\mathcal{M}.$ By the definition of
$\mathcal{M}$\ this means that $X_{-\gamma_{1}}X^{\alpha}\subset
(\mathcal{M)}_{-\gamma_{1}}=\overline{X_{-\gamma_{1}}I}^{\sigma}.$ It follows
that $I=J,$ so $J$ is $\mathcal{F}$-simple and%
\[
\mathcal{M=}\overline{\sum_{\gamma\in sp(\alpha),\gamma>0}X_{-\gamma
}J+X^{\alpha}([0,\infty))}^{\sigma}%
\]
Now let $\gamma\in sp(\alpha),\gamma\geqslant\gamma_{1}.$ Since $W=X,$\ it
follows that in particular
\[
X_{-\gamma}X^{\alpha}\subset(\mathcal{M)}_{-\gamma}=\overline{X_{-\gamma}%
J}^{\sigma}.
\]
By multiplying the previous relation to the left by $X_{\gamma}$\ we get%
\[
X_{\gamma}X_{-\gamma}X^{\alpha}=X_{\gamma}\overline{X_{-\gamma}J}^{\sigma
}\subset\overline{X_{\gamma}X_{-\gamma}J}^{\sigma}=J.
\]
Hence $\overline{X_{\gamma}X_{-\gamma}}^{\sigma}=J,$ so Condition (\textbf{S})
and therefore Condition 4.3. is satisfied.
\end{proof}

Next we will study the case when $0$\ is an accumulation point of $sp(\alpha)$.

Recall that all dynamical systems considered in this Section are supposed to
be non trivial, that is $sp(\alpha)\neq\left\{  0\right\}  $.

\bigskip

\textbf{4.12. Remark }\textit{Let (}$X,%
\mathbb{R}
,\alpha$\textit{) be an one-parameter }$\mathcal{F}$\textit{-dynamical system
such that }$X^{\alpha}([0,\infty))$\textit{ is a maximal }$\mathcal{F}%
$\textit{-closed subalgebra of }$X$\textit{ and let }$Y\subset X$ be as in
Lemma 4.10.\textit{\ If }%
\[
\overline{X^{\alpha}((-\infty,0))XX^{\alpha}((0,\infty))}^{\sigma}\subset
Y^{\alpha}.
\]
\textit{Then }%
\[
\overline{X^{\alpha}((0,\infty))XX^{\alpha}((-\infty,0))}^{\sigma}\subset
Y^{\alpha}.
\]
\textit{and conversely.}

\bigskip

\begin{proof}
\ First notice that $X^{\alpha}((0,\infty))^{2}=\left\{  0\right\}  .$ Indeed,
let $x,y\in X^{\alpha}((0,\infty))$ be such that $xy\neq0.$ Then $x^{\ast
}xy\neq0.$ But, according to the hypotheses, we have on the one hand
\[
x^{\ast}xy=x^{\ast}(xy)\in X^{\alpha}((-\infty,0))X^{\alpha}((0,\infty
))\subset X^{\alpha}.
\]
and on the other hand%
\[
x^{\ast}xy=(x^{\ast}x)y\in X^{\alpha}X^{\alpha}((0,\infty))\subset X^{\alpha
}((0,\infty)).
\]
Since $X^{\alpha}((0,\infty))\cap X^{\alpha}=\left\{  0\right\}  ,$ this shows
that $xy=0.$ Suppose that
\[
Z_{1}=\overline{X^{\alpha}((0,\infty))XX^{\alpha}((-\infty,0))}^{\sigma
}\nsubseteqq Y^{\alpha}.
\]
By the definition of $Z_{1}$\ it follows that $Z_{1}$\ is an $\alpha
$-invariant (hereditary) C*-subalgebra of $Y$\ and, by the above condition,
$sp(\alpha|_{Z_{1}})\neq\left\{  0\right\}  $, so $sp(\alpha|_{Z_{1}})$
contains negative numbers. Hence
\[
Z_{1}\nsubseteq X^{\alpha}([0,\infty)).
\]
\ Notice that
\[
X^{\alpha}([0,\infty))Z_{1}\subset Z_{1}.
\]
Since, obviously, $Z_{1}Z_{1}\subset Z_{1},$ and, as noticed above\ $Z_{1}%
\nsubseteq X^{\alpha}([0,\infty)),$\ it follows that the $\mathcal{F}$-closed
subalgebra $W$ of $X$ defined by
\[
W=\left\{  x\in X:xZ_{1}\subset Z_{1}\right\}  .
\]
strictly contains $X^{\alpha}([0,\infty))$ and therefore $W=X.$\ Thus $XZ_{1}$
$\subset Z_{1}$ and, as $Z_{1}$ is a C*-subalgebra, we also have
$Z_{1}X\subset Z_{1},$ so $XZ_{1}X\subset Z_{1}.$ But since $Y$\ is $\alpha
$-simple and $Z_{1}\subset Y,$\ it follows that$\ \overline{XZ_{1}X}^{\sigma
}=Y\subset Z_{1}$ so $Z_{1}=Y.$ On the other hand, since as shown above,
$X^{\alpha}((0,\infty))^{2}=\left\{  0\right\}  ,$ we have
\[
X^{\alpha}((0,\infty))Z_{1}=\left\{  0\right\}  .
\]
Thus
\[
X^{\alpha}((0,\infty))Y=\left\{  0\right\}  .
\]
Since $X^{\alpha}((0,\infty))=Y^{\alpha}((0,\infty)),$ it follows that
$X^{\alpha}((0,\infty))=\left\{  0\right\}  $ which is a contradiction with
our standing assumption $sp(\alpha)\neq\left\{  0\right\}  .$
\end{proof}

\bigskip

\textbf{4.13. Lemma }\textit{Suppose that\ (}$X,%
\mathbb{R}
,\alpha$\textit{) is a one-parameter }$\mathcal{F}$\textit{-dynamical system
such that }$sp(\alpha)$\textit{\ contains more than three points and
}$X^{\alpha}([0,\infty))$\textit{ is maximal.\ We have}%
\[
Z_{1}=\overline{X^{\alpha}((-\infty,0))XX^{\alpha}((0,\infty))}^{\sigma
}\nsubseteqq Y^{\alpha}.
\]
\textit{and, therefore, by Remark 4.12.}%
\[
Z_{2}=\overline{X^{\alpha}((0,\infty))XX^{\alpha}((-\infty,0))}^{\sigma
}\nsubseteqq Y^{\alpha}.
\]

\bigskip

\begin{proof}
Suppose to the contrary that
\[
\overline{X^{\alpha}((-\infty,0))XX^{\alpha}((0,\infty))}^{\sigma}\subset
Y^{\alpha}.
\]
so, in particular
\[
\overline{X^{\alpha}((-\infty,0))X^{\alpha}((0,\infty))}^{\sigma}\subset
Y^{\alpha}.
\]
Let $\gamma_{1},\gamma_{2}\in sp(\alpha),0<\gamma_{1}<\gamma_{2}$\ Then for
every $\epsilon>0,$ sufficiently small we have
\[
\overline{X^{\alpha}((-\gamma_{1}-\epsilon,-\gamma_{1}+\epsilon))X^{\alpha
}((\gamma_{1}-\epsilon,\gamma_{1}+\epsilon))}^{\sigma}\subset Y^{\alpha}.
\]
It follows that, for every $0<\epsilon<\frac{\gamma_{1}}{2},$ the set
\[
J_{\gamma_{1}}^{\epsilon}=\overline{X^{\alpha}((-\gamma_{1}-\epsilon
,-\gamma_{1}+\epsilon))X^{\alpha}((\gamma_{1}-\epsilon,\gamma_{1}+\epsilon
))}^{\sigma}%
\]
is an $\mathcal{F}$-closed ideal of $Y^{\alpha}.$ Let
\[
\mathcal{M}=\overline{J_{\gamma_{1}}^{\epsilon}X^{\alpha}((-\infty
,0))+X^{\alpha}([0,\infty))}^{\sigma}.
\]
Then $\mathcal{M}$\ is an $\mathcal{F}$-closed $\alpha$-invariant subspace of
$X.$ Denote
\[
W=\left\{  x\in X:x\mathcal{M}\subset\mathcal{M}\right\}  .
\]
Clearly, $W$\ is an $\mathcal{F}$-closed subalgebra of $X.$ Next, we will show
that $X^{\alpha}([0,\infty))\subset W.$ Clearly
\[
X^{\alpha}([0,\infty))X^{\alpha}([0,\infty))\subset X^{\alpha}([0,\infty
))\subset\mathcal{M}.
\]
so we have to prove that
\[
X^{\alpha}([0,\infty))J_{\gamma_{1}}^{\epsilon}X^{\alpha}((-\infty
,0))\subset\mathcal{M}.
\]
Indeed, by Lemma 2.11. ii)
\[
X^{\alpha}([0,\infty))=\cap_{\epsilon^{\prime}>0}X^{\alpha}((-\epsilon
^{\prime},\infty)).
\]
and by Lemma 2.11. iii)%
\[
X^{\alpha}((-\epsilon^{\prime},\infty))=\overline{X^{\alpha}((-\epsilon
^{\prime},\epsilon^{\prime}))+X^{\alpha}((\frac{\epsilon^{\prime}}{2}%
,\infty))}^{\sigma}.
\]
for every $\epsilon^{\prime}>0.$\ Hence, in order to prove that $X^{\alpha
}([0,\infty))J_{\gamma_{1}}^{\epsilon}X^{\alpha}((-\infty,0))\subset
\mathcal{M}$\ it is sufficient to prove that
\[
X^{\alpha}((-\epsilon,\epsilon))J_{\gamma_{1}}^{\epsilon}X^{\alpha}%
((-\infty,0))\subset\mathcal{M}\text{ and }X^{\alpha}((\frac{\epsilon}%
{2},\infty))J_{\gamma_{1}}^{\epsilon}X^{\alpha}((-\infty,0))\subset\mathcal{M}%
\]
for all $\epsilon>0$ sufficiently small. By Remark 4.12. we have
\[
X^{\alpha}((\frac{\epsilon}{2},\infty))J_{\gamma_{1}}^{\epsilon}X^{\alpha
}((-\infty,0))\subset X^{\alpha}X^{\alpha}((\gamma_{1}-\epsilon,\gamma
_{1}+\epsilon))X^{\alpha}((-\infty,0))\subset X^{\alpha}\subset\mathcal{M}.
\]
On the other hand%
\[
X^{\alpha}((-\epsilon,\epsilon))J_{\gamma_{1}}^{\epsilon}X^{\alpha}%
((-\infty,0))\mathcal{=}%
\]%
\[
X^{\alpha}((-\epsilon,\epsilon))\overline{X^{\alpha}((-\gamma_{1}%
-\epsilon,-\gamma_{1}+\epsilon))X^{\alpha}((\gamma_{1}-\epsilon,\gamma
_{1}+\epsilon))}^{\sigma}J_{\gamma_{1}}^{\epsilon}X^{\alpha}((-\infty
,0))\subset
\]%
\[
\subset X^{\alpha}((-\gamma_{1}-2\epsilon,-\gamma_{1}+2\epsilon))X^{\alpha
}((\gamma_{1}-\epsilon,\gamma_{1}+\epsilon))J_{\gamma_{1}}^{\epsilon}%
X^{\alpha}((-\infty,0))\subset
\]%
\[
\subset X^{\alpha}J_{\gamma_{1}}^{\epsilon}X^{\alpha}((-\infty,0))\subset
J_{\gamma_{1}}^{\epsilon}X^{\alpha}((-\infty,0))\subset\mathcal{M}.
\]
We have thus proven that $X^{\alpha}([0,\infty))\subset W.$ Now we will show
that this inclusion is strict. Indeed,\ let $0\neq x\in X^{\alpha}%
((-\gamma_{1}-\epsilon,-\gamma_{1}+\epsilon))$ and denote $y=xx^{\ast}x\in
J_{\gamma_{1}}^{\epsilon}X^{\alpha}((-\gamma_{1}-\epsilon,-\gamma_{1}%
+\epsilon))\subset X^{\alpha}((-\gamma_{1}-\epsilon,-\gamma_{1}+\epsilon)).$
Then, $y\neq0$\ and $y\notin X^{\alpha}([0,\infty)).$ We have
\[
yJ_{\gamma_{1}}^{\epsilon}X^{\alpha}((-\infty,0))=xx^{\ast}xJ_{\gamma_{1}%
}^{\epsilon}X^{\alpha}((-\infty,0))\subset J_{\gamma_{1}}^{\epsilon}X^{\alpha
}((-\infty,0))^{2}=\left\{  0\right\}  .
\]
We show next that $yX^{\alpha}([0,\infty))\subset\mathcal{M}.$ As noticed
above, we have
\[
X^{\alpha}([0,\infty))=\cap_{\epsilon^{^{\prime}}>0}X^{\alpha}((-\epsilon
^{^{\prime}},\infty)).
\]
and
\[
X^{\alpha}((-\epsilon^{^{\prime}},\infty))=\overline{X^{\alpha}((-\epsilon
^{^{\prime}},\epsilon^{^{\prime}}))+X^{\alpha}((\frac{\epsilon^{^{\prime}}}%
{2},\infty))}^{\sigma}.
\]
Clearly, if $0<\epsilon^{^{\prime}}<\frac{\gamma_{1}}{2}$ we have%
\[
yX^{\alpha}((-\epsilon^{^{\prime}},\epsilon^{^{\prime}}))\subset J_{\gamma
_{1}}^{\epsilon}X^{\alpha}((-\gamma_{1}-\epsilon-\epsilon^{^{\prime}}%
,-\gamma_{1}+\epsilon+\epsilon^{^{\prime}}))\subset J_{\gamma_{1}}^{\epsilon
}X^{\alpha}((-\infty,0))\subset\mathcal{M}.
\]
and%
\[
yX^{\alpha}((\frac{\epsilon^{^{\prime}}}{2},\infty))=(xx^{\ast})xX^{\alpha
}((\frac{\epsilon^{^{\prime}}}{2},\infty))\subset J_{\gamma_{1}}^{\epsilon
}X^{\alpha}((-\infty,0))X^{\alpha}((0,\infty))\subset X^{\alpha}%
\subset\mathcal{M}.
\]
So $X^{\alpha}([0,\infty))\varsubsetneq W.$ Since $X^{\alpha}([0,\infty))$ is
maximal, it follows that $W=X.$ In order to show that the assumption
\[
\overline{X^{\alpha}((-\infty,0))ZX^{\alpha}((0,\infty))}^{\sigma}\subset
Y^{\alpha}%
\]
leads to a contradiction,\ we will use the hypothesis that there is
$\gamma_{2}\in sp(\alpha)$ such that $0<\gamma_{1}<\gamma_{2}.$\ Let
$0<\epsilon<\min\left\{  \frac{\gamma_{1}}{2},\frac{\gamma_{2}-\gamma_{1}}%
{2}\right\}  .$\ Since $W=X$\ \ we have%
\[
X^{\alpha}((-\gamma_{2}-\epsilon,-\gamma_{2}+\epsilon))\mathcal{M\subset M}.
\]
In particular%
\[
X^{\alpha}((-\gamma_{2}-\epsilon,-\gamma_{2}+\epsilon))X^{\alpha}%
\subset\mathcal{M}.
\]
so\
\[
\left\{  0\right\}  \neq X^{\alpha}((-\gamma_{2}-\epsilon,-\gamma_{2}%
+\epsilon))X^{\alpha}((\gamma_{2}-\epsilon,\gamma_{2}+\epsilon))X^{\alpha
}((-\gamma_{2}-\epsilon,-\gamma_{2}+\epsilon))\subset\mathcal{M}.
\]
If we denote $J_{\gamma_{2}}^{\epsilon}=\overline{X^{\alpha}((-\gamma
_{2}-\epsilon,-\gamma_{2}+\epsilon))X^{\alpha}((\gamma_{2}-\epsilon,\gamma
_{2}+\epsilon))}^{\sigma}$\ then $J_{\gamma_{2}}^{\epsilon}$ is an ideal of
$X^{\alpha}$ and%
\[
J_{\gamma_{2}}^{\epsilon}X^{\alpha}((-\gamma_{2}-\epsilon,-\gamma_{2}%
+\epsilon))\subset\mathcal{M}%
\]
At the same time, since $\mathcal{M}$\ is an $\mathcal{F}$-closed $\alpha
$-invariant\ subspace of $X,$ and $J_{\gamma_{2}}^{\epsilon}X^{\alpha
}((-\gamma_{2}-\epsilon,-\gamma_{2}+\epsilon))\subset X^{\alpha}%
((-\infty,0)),$\ it follows that%
\[
J_{\gamma_{2}}^{\epsilon}X^{\alpha}((-\gamma_{2}-\epsilon,-\gamma_{2}%
+\epsilon))\subset\mathcal{M}^{\alpha}((-\infty.0))=\overline{J_{\gamma_{1}%
}^{\epsilon}X^{\alpha}((-\infty,0))}^{\sigma}.
\]
Multiplying to the right the above equality by $X^{\alpha}((\gamma
_{2}-\epsilon,\gamma_{2}+\epsilon))$ we get%
\[
J_{\gamma_{2}}^{\epsilon}\subset J_{\gamma_{1}}^{\epsilon}%
\]
so, since $J_{\gamma_{2}}^{\epsilon}$ and $J_{\gamma_{1}}^{\epsilon}$ are
ideals of $X^{\alpha},$\ it follows that $J_{\gamma_{2}}^{\epsilon}%
J_{\gamma_{1}}^{\epsilon}=J_{\gamma_{2}}^{\epsilon}\neq\left\{  0\right\}  .$
Using the definition of the ideals $J_{\gamma_{2}}^{\epsilon}$ and
$J_{\gamma_{1}}^{\epsilon}$ we get%
\[
X^{\alpha}((\gamma_{2}-\epsilon,\gamma_{2}+\epsilon))X^{\alpha}((-\gamma
_{1}-\epsilon,-\gamma_{1}+\epsilon))\neq\left\{  0\right\}  .
\]
and at the same time%
\[
X^{\alpha}((\gamma_{2}-\epsilon,\gamma_{2}+\epsilon))X^{\alpha}((-\gamma
_{1}-\epsilon,-\gamma_{1}+\epsilon))\subset X^{\alpha}%
\]
This is a contradiction since $X^{\alpha}((\gamma_{2}-\epsilon,\gamma
_{2}+\epsilon))X^{\alpha}((-\gamma_{1}-\epsilon,-\gamma_{1}+\epsilon))\subset
X^{\alpha}((\gamma_{2}-\gamma_{1}-2\epsilon,\gamma_{2}-\gamma_{1}%
+2\epsilon))\subset X^{\alpha}((0,\infty))$\ and the proof is completed.\ 
\end{proof}

\bigskip

\textbf{4.14. Lemma }\textit{Let (}$X,%
\mathbb{R}
,\alpha$\textit{) be a one-parameter }$\mathcal{F}$\textit{-dynamical system
such that\ }$sp(\alpha)$\textit{\ contains more than three points and let
}$Z\in H_{\sigma}^{\alpha}.$\textit{\ If }$X^{\alpha}([0,\infty))$\textit{\ is
a maximal }$\mathcal{F}$\textit{-closed subalgebra of }$X,$\textit{ then, }%
\[
\overline{X^{\alpha}((-\infty,-\gamma))ZX^{\alpha}((\gamma,\infty))}^{\sigma
}=Y.
\]
\textit{and }%
\[
\overline{X^{\alpha}((\gamma,\infty))ZX^{\alpha}((-\infty,-\gamma))}^{\sigma
}=Y.
\]
\textit{for every }$\gamma\in%
\mathbb{R}
,\gamma\geqslant0$\textit{ such that }$sp(\alpha)\cap(\gamma,\infty
)\neq\emptyset$\textit{ where }$Y$\textit{\ is the }$\alpha$\textit{-simple
ideal of }$X$\textit{\ whose existence was established in Lemma 4.10.}

\bigskip

\begin{proof}
Let $\gamma\in%
\mathbb{R}
,\gamma\geqslant0$\textit{ be such that }$sp(\alpha)\cap(\gamma,\infty
)\neq\emptyset,$ so $X^{\alpha}((\gamma,\infty))\neq\left\{  0\right\}  .$
Consider first the case when $Z=X.$\ Let $Z_{1}=\overline{X^{\alpha}%
((\gamma,\infty))XX^{\alpha}((-\infty,\gamma))}^{\sigma}.$ Since $Y$\ is
$\alpha$-simple and $Z_{1}$\ is $\alpha$-invariant, we have that
$\overline{YZ_{1}Y}^{\sigma}=Y,$\ so, if $Z_{1}\subsetneq Y,$ we have
$\overline{YZ_{1}}^{\sigma}\nsubseteq Z_{1}.$\ Therefore, if $Z_{1}\subsetneq
Y,$\ then, there exist $y\in Y$\ and $z\in Z_{1}$\ such that $yz\notin Z_{1}.$
Denote%
\[
W=\left\{  x\in X:xZ_{1}\subset Z_{1}\right\}  .
\]
Then, $W$ is a $\mathcal{F}-$ closed subalgebra of $X$\ and, by the above
argument, if $Z_{1}\subsetneq Y,$we have that $W\neq X.$\ Let us prove that
\[
X^{\alpha}([0,\infty))\subset W.
\]
By Lemma 2.11. i), $X^{\alpha}((\gamma,\infty))=\overline{\sum_{n\in%
\mathbb{N}
}X^{\alpha}([\gamma+\frac{1}{n},\infty))}^{\sigma}.$ Applying Lemma 2.11. v),
for every $n\in%
\mathbb{N}
,$\ we have
\[
X^{\alpha}([0,\infty))X^{\alpha}([\gamma+\frac{1}{n},\infty))\subset
X^{\alpha}([\gamma+\frac{1}{n},\infty)).
\]
Therefore,%
\[
X^{\alpha}([0,\infty))\sum_{n\in%
\mathbb{N}
}X^{\alpha}([\gamma+\frac{1}{n},\infty))\subset\sum_{n\in%
\mathbb{N}
}X^{\alpha}([\gamma+\frac{1}{n},\infty)).
\]
So
\[
X^{\alpha}([0,\infty))X^{\alpha}((\gamma,\infty))\subset X^{\alpha}%
((\gamma,\infty)).
\]
and therefore%
\[
X^{\alpha}([0,\infty))Z_{1}\subset Z_{1}.
\]
and hence
\[
X^{\alpha}([0,\infty))\subset W.
\]
On the other hand, clearly, $Z_{1}\subset W$\ and, since $sp(\alpha)$ contains
at least five points$,$ by the previous Lemma 4.13. it follows that
$Z_{1}\nsubseteq Y^{\alpha}$, so $sp(\alpha|_{Z_{1}})\neq\left\{  0\right\}
.$ Since $Z_{1}$\ is an $\alpha$-invariant (hereditary) C*-subalgebra of
\ there exists $\gamma\in sp(\alpha|_{Z_{1}}),\gamma<0,$ so $Z_{1}-X^{\alpha
}([0,\infty))\neq\emptyset.$\ Therefore, if $z\in Z_{1}-X^{\alpha}%
([0,\infty)),$\ then $z\in W-$ $X^{\alpha}([0,\infty)).$ We have thus proved
that if $Z\subsetneq Y,$\ then there exists an $\mathcal{F}$-closed subalgebra
$W$\ of $X$\ such that%
\[
X^{\alpha}([0,\infty))\subsetneq W\subsetneq X.
\]
which contradicts the maximality of $X^{\alpha}([0,\infty)).$ Hence $Z_{1}=Y.$
To prove that
\[
\overline{X^{\alpha}((-\infty,-\gamma))XX^{\alpha}((\gamma,\infty))}^{\sigma
}=Y.
\]
let $Z_{1}=\overline{X^{\alpha}((-\infty,-\gamma))XX^{\alpha}((\gamma
,\infty))}^{\sigma}.$ If we denote%
\[
W=\left\{  x\in X:Z_{1}x\subset Z_{1}\right\}  .
\]
then, by similar arguments it can be proven that $Z_{1}=Y.$ Now let
$Z\in\mathcal{H}_{\sigma}^{\alpha}(Y)$ be arbitrary and denote%
\[
Z_{1}=\overline{X^{\alpha}((\gamma,\infty))ZX^{\alpha}((-\infty,-\gamma
))}^{\sigma}.
\]
Then, notice that $Z_{1}\neq\left\{  0\right\}  .$ Indeed if $Z_{1}=\left\{
0\right\}  ,$ then $X^{\alpha}((0,\infty))Z=\left\{  0\right\}  .$ By the
first part of the proof, it follows that $YZ=\left\{  0\right\}  ,$ so
$Z=\left\{  0\right\}  ,$ which contradicts the fact that $Z\in\mathcal{H}%
_{\sigma}^{\alpha}(Y).$ We notice also that $Z_{1}\nsubseteqq Y^{\alpha}.$
Indeed, suppose to the contrary that $Z_{1}\subset Y^{\alpha}.$ Then, since
$X^{\alpha}((0,\infty))Z_{1}\subset Z_{1}\subset Y^{\alpha},$ and, on the
other hand, $X^{\alpha}((0,\infty))Z_{1}\subset Y^{\alpha}((0,\infty
)),$\ it\ follows that $X^{\alpha}((0,\infty))Z_{1}=\left\{  0\right\}  .$
Therefore
\[
X^{\alpha}((-\infty,0))XX^{\alpha}((0,\infty))Z_{1}=\left\{  0\right\}  .
\]
so $YZ_{1}=\left\{  0\right\}  $ and thus $Z_{1}=\left\{  0\right\}  $,
contradiction. Using these two facts about $Z_{1}$ for arbitrary
$Z\in\mathcal{H}_{\sigma}^{\alpha}(Y)$\ and the arguments above for the
particular case $Z=X$ the proof is completed.
\end{proof}

\bigskip

The\bigskip next lemma is probably known, but I include its proof below.

\textbf{4.15. Lemma }\textit{Let }$\mathit{Y}$\textit{ be a C*-algebra and
}$D\subset Y$ \textit{a hereditary} C*-\textit{subalgebra of} $Y$\textit{. Let
}$\left\{  A_{i}\right\}  _{i\in I}$\textit{ be a collection of subsets of
}$Y$\textit{ such that }$A_{i}^{\ast}YA_{i}\subset D$\textit{ for all }$i\in
I.$ \textit{Then, }$(\sum A_{i})^{\ast}Y(\sum A_{i})\subset D$ (Here, the
symbol $\ \sum A_{i}\ $denotes the set of all finite sums of $A_{i}^{\prime}%
s$)$.$

\bigskip

\begin{proof}
Let $i,j\in I,i\neq j$ be arbitrary$.$ If we prove that $A_{i}^{\ast}%
YA_{j}\subset D,$ it will follow that\textit{ } $(\sum A_{i})^{\ast}Y(\sum
A_{i})\subset D.$ Let $a\in A_{i},y\in Y$ and $b\in A_{j}.$ Then, using the
well known and easy to prove inequality%
\[
(a^{\ast}y^{\ast}\pm b^{\ast}y^{\ast})(ya\pm yb)\leq2(a^{\ast}y^{\ast
}ya+b^{\ast}y^{\ast}yb)
\]
and the hypotheses that $a^{\ast}y^{\ast}ya\in D,b^{\ast}y^{\ast}yb\in D$ and
that $D$\ is a hereditary C*-subalgebra of $Y,$ it follows that%
\[
a^{\ast}y^{\ast}yb,\pm b^{\ast}y^{\ast}ya\in D
\]
so $A_{i}^{\ast}Y^{+}A_{j}\subset D$ (where $Y^{+}$\ denotes the set of all
non negative elements of $Y)$\ and therefore $A_{i}^{\ast}YA_{j}\subset D.$\ 
\end{proof}

\bigskip

\textbf{4.16. Lemma }\textit{Let }$(X,%
\mathbb{R}
,\alpha)$\textit{\ be an }$\mathcal{F}$\textit{-dynamical system. Suppose that
}$0$\textit{\ is an accumulation point of }$sp(\alpha)$\textit{ and
}$X^{\alpha}([0,\infty))$\textit{\ is a maximal }$\mathcal{F}$\textit{-closed
subalgebra of }$X.$\textit{\ Then Condition 4.3. b) holds.}

\bigskip

\textbf{Proof. \ }Let $Y$\ be the $\alpha$-simple ideal of $X$\ from Lemma
4.10. We will prove this lemma in three steps:

\bigskip

\textbf{Step 1}. Let $\gamma_{0}\in sp(\alpha)$ and $\epsilon>0.$ Then%
\[
S_{\gamma_{0}}=\overline{X^{\alpha}((-\gamma_{0}-\epsilon,-\gamma_{0}%
+\epsilon))XX^{\alpha}((\gamma_{0}-\epsilon,\gamma_{0}+\epsilon))}^{\sigma
}=Y.
\]
and%
\[
T_{\gamma_{0}}=\overline{X^{\alpha}((\gamma_{0}-\epsilon,\gamma_{0}%
+\epsilon))XX^{\alpha}((-\gamma_{0}-\epsilon,-\gamma_{0}+\epsilon))}^{\sigma
}=Y.
\]
where $S_{\gamma_{0}}$ and $T_{\gamma_{0}}$ are notations for the
corresponding $\mathcal{F}$-closed hereditary subalgebras of $Y.$ We will
prove the first equality and then show how to obtain the second one.
Obviously, it is sufficient to prove the equality for $0<\epsilon<\frac
{\gamma_{0}}{2}.$\ Denote%
\[
\mathcal{M}=\overline{\overline{XS_{\gamma_{0}}}^{\sigma}+X^{\alpha}%
((\gamma_{0}-\frac{\epsilon}{2},\infty))}^{\sigma}%
\]
and%
\[
W=\left\{  x\in X:x\mathcal{M}\subset\mathcal{M}\right\}
\]
Then, $\mathcal{M}$\ is an $\mathcal{F}$-closed subspace of $X$\ and $W$\ is
an $\mathcal{F}$-closed subalgebra of $X.$ Using Lemma 2.11. i) and v) as in
the proof of Lemma 4.14. it follows immediately that $X^{\alpha}%
([0,\infty))\subset W.$ Next, we will show that $X^{\alpha}([0,\infty))$\ is a
proper subset of $W.$\ Indeed, since $0$\ is a point of accumulation of
$sp(\alpha),$\ there exists $0\neq x\in X^{\alpha}((0,\frac{\epsilon}%
{4}))=Y^{\alpha}((0,\frac{\epsilon}{4})).$ We will show that $x^{\ast}\in W.$
Since clearly $x^{\ast}\notin X^{\alpha}([0,\infty))$\ the above claim will be
proven. By Lemma 2.11. iii) we have%
\[
X^{\alpha}((\gamma_{0}-\frac{\epsilon}{2},\infty))=\overline{X^{\alpha
}((\gamma_{0}-\frac{\epsilon}{2},\gamma_{0}+\epsilon))+X^{\alpha}((\gamma
_{0}+\frac{3\epsilon}{4},\infty))}^{\sigma}.
\]
Taking into account that $Y$\ is an $\alpha$-simple ideal of $X$\ it follows
that%
\[
x^{\ast}X^{\alpha}((\gamma_{0}-\frac{\epsilon}{2},\gamma_{0}+\epsilon))\subset
Y^{\alpha}((\gamma_{0}-\frac{3\epsilon}{4},\gamma_{0}+\epsilon))\subset
Y^{\alpha}((\gamma_{0}-\epsilon,\gamma_{0}+\epsilon))\subset
\]%
\[
\subset\overline{YY^{\alpha}((\gamma_{0}-\epsilon,\gamma_{0}+\epsilon
))}^{\sigma}=\overline{YY^{\alpha}((-\gamma_{0}-\epsilon,-\gamma_{0}%
+\epsilon))YY^{\alpha}((\gamma_{0}-\epsilon,\gamma_{0}+\epsilon))}^{\sigma
}\subset
\]%
\[
\subset\overline{YS_{\gamma_{0}}}^{\sigma}.
\]
and%
\[
x^{\ast}X^{\alpha}((\gamma_{0}+\frac{3\epsilon}{4},\infty))\subset X^{\alpha
}((\gamma_{0}-\frac{\epsilon}{2},\infty)).
\]
Therefore $x^{\ast}\in W-X^{\alpha}([0,\infty)),$ so $X^{\alpha}%
([0,\infty))\varsubsetneq W.$ Since $X^{\alpha}([0,\infty))$ is a maximal
$\mathcal{F}$-closed subalgebra of $X$\ it follows that $W=X.$ We will prove
that $S_{\gamma_{0}}=Y$ by contrapositive. Suppose that $S_{\gamma_{0}%
}\varsubsetneq Y.$ We will produce an element $x\neq0$ such that $x^{\ast}\in
X-W.$ Since $(\gamma_{0}-\frac{\epsilon}{2},\infty)=\cup_{\gamma>\gamma
_{0}-\frac{\epsilon}{2}}(\gamma,\gamma+\gamma_{0}-\frac{\epsilon}{2}),$
from\ Lemma 2.11. iii) it follows that%
\[
X^{\alpha}((\gamma_{0}-\frac{\epsilon}{2},\infty))=\overline{\sum
_{\gamma>\gamma_{0}-\frac{\epsilon}{2}}X^{\alpha}((\gamma,\gamma+\gamma
_{0}-\frac{\epsilon}{2}))}^{\sigma}.
\]
By Lemma 4.14.
\[
\overline{X^{\alpha}((\gamma_{0}-\frac{\epsilon}{2},\infty))^{\ast}YX^{\alpha
}((-\infty,-\gamma_{0}+\frac{\epsilon}{2}))}^{\sigma}=Y.
\]
Hence%
\[
\overline{(\sum_{\gamma>\gamma_{0}-\frac{\epsilon}{2}}X^{\alpha}%
((\gamma,\gamma+\gamma_{0}-\frac{\epsilon}{2})))^{\ast}Y(\sum_{\gamma
>\gamma_{0}-\frac{\epsilon}{2}}X^{\alpha}((\gamma,\gamma+\gamma_{0}%
-\frac{\epsilon}{2})))}^{\sigma}=Y.
\]
Since we are assuming that $S_{\gamma_{0}}\varsubsetneq Y,$ Lemma 4.15.
implies that there exists $\gamma>\gamma_{0}-\frac{\epsilon}{2}$\ such that%
\[
X^{\alpha}((\gamma,\gamma+\gamma_{0}-\frac{\epsilon}{2}))^{\ast}YX^{\alpha
}((\gamma,\gamma+\gamma_{0}-\frac{\epsilon}{2}))^{\ast}\nsubseteq
S_{\gamma_{0}}.
\]
Since $X^{\alpha}((\gamma,\gamma+\gamma_{0}-\frac{\epsilon}{2}))$ is a
subspace of $X,$ there exists $x\in X^{\alpha}((\gamma,\gamma+\gamma_{0}%
-\frac{\epsilon}{2}))$ such that $x^{\ast}x\notin S_{\gamma_{0}}.$ Then
$x^{\ast}x\notin YS_{\gamma_{0}}$\ since otherwise ($x^{\ast}x)^{2}\in
S_{\gamma_{0}}YS_{\gamma_{0}}\subset S_{\gamma_{0}},$ so $x^{\ast}x\in
S_{\gamma_{0}}.$ On the other hand
\[
x^{\ast}x\in Y^{\alpha}((-\gamma_{0}+\frac{\epsilon}{2},\gamma_{0}%
-\frac{\epsilon}{2})).
\]
so $x^{\ast}x\notin X^{\alpha}((\gamma_{0}-\frac{\epsilon}{2},\infty)).$
Notice that both $\overline{YS_{\gamma_{0}}}^{\sigma}$\ and $X^{\alpha
}((\gamma_{0}-\frac{\epsilon}{2},\infty))$\ are $\alpha$-invariant subspaces,
so if $z\in$\ $\overline{YS_{\gamma_{0}}}^{\sigma}$\ (respectively $y\in
X^{\alpha}((\gamma_{0}-\frac{\epsilon}{2},\infty))$) and $f\in L^{1}(%
\mathbb{R}
)$\ we have $\alpha_{f}(z)\in$ $\overline{YS_{\gamma_{0}}}^{\sigma}$
(respectively $\alpha_{f}(y)\in X^{\alpha}((\gamma_{0}-\frac{\epsilon}%
{2},\infty))$). Next, we will prove that $x^{\ast}x\notin\mathcal{M},$\ so
\ \ $x^{\ast}\notin W.$ Indeed if $x^{\ast}x\in\mathcal{M},$\ then there
exists a net $\left\{  z_{\iota}+y_{\iota}:z_{\iota}\in YS_{\gamma_{0}%
},y_{\iota}\in X^{\alpha}((\gamma_{0}-\frac{\epsilon}{2},\infty))\right\}
$\ such that $x^{\ast}x=\mathcal{F}$-$lim(z_{\iota}+y_{\iota}).$ By [20,
Theorem 2.6.2.] there exists a function $f\in L^{1}(%
\mathbb{R}
)$\ such that $\widehat{f}=1$ on an open set containing $sp(x^{\ast}x)$ whose
closure is included in $(-\gamma_{0}+\frac{\epsilon}{2},\gamma_{0}%
-\frac{\epsilon}{2})$\ and $\widehat{f}=0$ outside $(-\gamma_{0}%
+\frac{\epsilon}{2},\gamma_{0}-\frac{\epsilon}{2}).$ Hence, since $\alpha_{f}%
$\ is $\mathcal{F}$-continuous, we have
\[
\alpha_{f}(x^{\ast}x)=x^{\ast}x=\mathcal{F}-\lim\alpha_{f}(z_{\iota}+y_{\iota
})=\mathcal{F}-\lim\alpha_{f}(z_{\iota})\in\overline{YS_{\gamma_{0}}}^{\sigma
}.
\]
contradiction. Hence if $S_{\gamma_{0}}\varsubsetneq Y,$ then $W\varsubsetneq
X$ so $X^{\alpha}([0,\infty))$\textit{\ }is not maximal\textit{.} To prove
that $T_{\gamma_{0}}=Y,$ replace the above $\mathcal{M}$\ and $W$\ by
\[
\mathcal{M}=\overline{\overline{T_{\gamma_{0}}Y}^{\sigma}+X^{\alpha}%
((\gamma_{0}-\frac{\epsilon}{2},\infty))}^{\sigma}.
\]%
\[
W=\left\{  x\in X:\mathcal{M}x\subset\mathcal{M}\right\}  .
\]
and use similar arguments.

\bigskip

\textbf{Step 2.} Let $\gamma_{0}\in sp(\alpha)$ and $0<\epsilon<\frac
{\gamma_{0}}{2}.\ $If $Z\in\mathcal{H}_{\sigma}^{\alpha}(Y),$ then%
\[
S_{\gamma_{0}}^{\epsilon}(Z)=\overline{X^{\alpha}((-\gamma_{0}-\epsilon
,-\gamma_{0}+\epsilon))ZX^{\alpha}((\gamma_{0}-\epsilon,\gamma_{0}+\epsilon
))}^{\sigma}=Y.
\]
and%
\[
T_{\gamma_{0}}^{\epsilon}(Z)=\overline{X^{\alpha}((\gamma_{0}-\epsilon
,\gamma_{0}+\epsilon))ZX^{\alpha}((-\gamma_{0}-\epsilon,-\gamma_{0}%
+\epsilon))}^{\sigma}=Y.
\]
We will prove only the first equality and then describe how to obtain the
second one. As in Step 1. denote%
\[
\mathcal{M}=\overline{\overline{XS_{\gamma_{0}}(Z)}^{\sigma}+X^{\alpha
}((\gamma_{0}-\frac{\epsilon}{2},\infty))}^{\sigma}.
\]
and%
\[
W=\left\{  x\in X:x\mathcal{M}\subset\mathcal{M}\right\}  .
\]
Then $\mathcal{M}$\ is an $\mathcal{F}$-closed subspace of $X$\ and $W$\ is an
$\mathcal{F}$-closed subalgebra of $X.$ Clearly $X^{\alpha}([0,\infty))\subset
W.$ The fact that this inclusion is strict is the only difference between the
case of general $Z\in\mathcal{H}_{\sigma}^{\alpha}(Y)$ and the case when $Z=X$
considered in Step 1. Since $0$\ is an accumulation point of $sp(\alpha),$
there exists $\gamma_{1}\in sp(\alpha),\gamma_{1}>0$ such that $\gamma
_{1}<\frac{\epsilon}{8}.$ By Step 1. we have%
\[
\overline{X^{\alpha}((-\gamma_{1}-\epsilon^{\prime},-\gamma_{1}+\epsilon
^{\prime}))YX^{\alpha}((\gamma_{1}-\epsilon^{\prime},\gamma_{1}+\epsilon
^{\prime}))}^{\sigma}=Y
\]
and%
\[
\overline{X^{\alpha}((\gamma_{1}-\epsilon^{\prime},\gamma_{1}+\epsilon
^{\prime}))YX^{\alpha}((-\gamma_{1}-\epsilon^{\prime},-\gamma_{1}%
+\epsilon^{\prime}))}^{\sigma}=Y
\]
for every $\epsilon^{\prime}>0.$\ Hence $S_{\gamma_{1}}^{\epsilon}%
(Z)\neq\left\{  0\right\}  ,$ since, otherwise $YZY=\left\{  0\right\}  ,$
hence $Z=\left\{  0\right\}  $ which is not possible since $Z\in
\mathcal{H}_{\sigma}^{\alpha}(Y),$ so $Z\neq\left\{  0\right\}  $ by the
definition of $\mathcal{H}_{\sigma}^{\alpha}(Y).$ In addition, for every
$\delta>0,$ and $\epsilon^{\prime}>0$\ we have%
\[
\overline{X^{\alpha}((-\gamma_{0}-\epsilon^{\prime},-\gamma_{0}+\epsilon
^{\prime}))Z^{\alpha}((-\delta,\delta))X^{\alpha}((\gamma_{0}-\epsilon
^{\prime},\gamma_{0}+\epsilon^{\prime}))}^{\sigma}\neq\left\{  0\right\}
\]
for the same reason and the fact that $0\in sp(\alpha|_{Z})$. Since
$\gamma_{1}<\frac{\epsilon}{8}$ there exists $n\in%
\mathbb{N}
,n>4\ $such that $\frac{\epsilon}{2^{n}}\leq\gamma_{1}<\frac{\epsilon}%
{2^{n-1}}.$ Then, $\left\{  0\right\}  \neq Z^{\alpha}(-\frac{\epsilon
}{2^{n+1}},\frac{\epsilon}{2^{n+1}})X^{\alpha}(\gamma_{1}-\frac{\epsilon
}{2^{n+1}},\gamma_{1}+\frac{\epsilon}{2^{n+1}})\subset X^{\alpha}%
(0,\frac{\epsilon}{4}).$ Let $x\in Z^{\alpha}(-\frac{\epsilon}{2^{n+1}}%
,\frac{\epsilon}{2^{n+1}})X^{\alpha}(\gamma_{1}-\frac{\epsilon}{2^{n+1}%
},\gamma_{1}+\frac{\epsilon}{2^{n+1}}),x\neq0.$ Clearly, $x^{\ast}\notin
X^{\alpha}([0,\infty)).$ The proof that $x^{\ast}\in W\ $is very similar with
the corresponding proof in Step 1. and we will omit it. The rest of the proof
is a verbatim repetition of the arguments in Step 1.

\bigskip

\textbf{Step 3.} Proof of Lemma. Applying Remark 4.2. v) to $(Z,%
\mathbb{R}
,\alpha)$\ and $\epsilon>0,$\ it follows that%

\[
\overline{Z^{\alpha}((-\epsilon,\epsilon))ZZ^{\alpha}((-\epsilon,\epsilon
))}^{\sigma}=Z\text{ \ \ \ \ \ \ \ \ (1)}%
\]
and therefore, since $Z$\ is a hereditary C*-subalgebra of $Y$%
\[
\overline{Z^{\alpha}((-\epsilon,\epsilon))YZ^{\alpha}((-\epsilon,\epsilon
))}^{\sigma}=Z\text{ \ \ \ \ \ \ \ \ \ \ \ \ \ (2)}%
\]
Let $\gamma_{0}\in sp(\alpha)-\left\{  0\right\}  $ and $\epsilon>0$ as in
Step 2. By replacing $Y$ in equation (2) above by $\overline{X^{\alpha
}((-\gamma_{0}-\epsilon,-\gamma_{0}+\epsilon))ZX^{\alpha}((\gamma_{0}%
-\epsilon,\gamma_{0}+\epsilon))}^{\sigma}$ as in Step 2 and then replacing $Z$
by $\overline{Z^{\alpha}((-\epsilon,\epsilon))ZZ^{\alpha}((-\epsilon
,\epsilon))}^{\sigma}$ as in relation (1) above, we get \ \ \ \
\[
\overline{Z^{\alpha}((-\gamma_{0}-3\epsilon,-\gamma_{0}+3\epsilon))ZZ^{\alpha
}((\gamma_{0}-3\epsilon,\gamma_{0}+3\epsilon))}^{\sigma}=Z
\]
Therefore%
\[
\overline{Z^{\alpha}((-\gamma_{0}-\delta,-\gamma_{0}+\delta))ZZ^{\alpha
}((\gamma_{0}-\delta,\gamma_{0}+\delta))}^{\sigma}=Z
\]
for all $\delta>0$ and consequently $sp(\alpha|_{Y})=\widetilde{\Gamma
}_{\mathcal{F}}(\alpha|_{Y}).\blacksquare$

\bigskip

\textbf{4.17. Corollary \ }\textit{Let (}$X,%
\mathbb{R}
,\alpha$\textit{) be a one-parameter }$\mathcal{F}$\textit{-dynamical system.
and such that }$X^{\alpha}([0,\infty))$\textit{\ is a maximal }$\mathcal{F}%
$\textit{-closed subalgebra of }$X$\textit{. Then Condition 4.3. is
satisfied.}

\bigskip

\begin{proof}
Follows from Lemmas 4.11. and 4.16.
\end{proof}

\bigskip

Now we can state the main result of this paper. This result contains and
improves on all the previous results about maximality of the algebra of
analytic elements associated with a C*- or W*-one-parameter dynamical system.
Moreover it also answers the maximality question for multiplier one-parameter
dynamical systems.

\bigskip

\textbf{4.18. Theorem }\textit{Let (}$X,%
\mathbb{R}
,\alpha$\textit{) be a one-parameter }$\mathcal{F}$\textit{-dynamical system.
Then }$X^{\alpha}([0,\infty))$\textit{\ is a maximal }$\mathcal{F}%
$\textit{-closed subalgebra of }$X$\textit{ if and only if the Spectral
Condition 4.3. is satisfied.}

\bigskip

\begin{proof}
Follows from Theorem 4.9. and Corollary 4.17.
\end{proof}

\bigskip

\textbf{Acknowledgment. }I would like to thank the referee for valuable
comments which helped to simplify some proofs and improve the paper.

\bigskip\ 

\begin{center}
\textbf{References}
\end{center}

\bigskip\lbrack1] W. B. ARVESON, On groups of automorphisms of operator
algebras, \textit{J. Funct. Anal. 15 (1974), 217-243.}

\bigskip

[2] W. B. ARVESON, The harmonic analysis of automorphism groups,
\textit{Operator Algebras and Applications, Part I, Kingston, Ontario 1980,
Proc. Symposia Pure Math., Vol. 38, AMS, Providence, RI 1982, 199-269.}

\bigskip

[3] A. CONNES, Une classification des facteurs de type III, \textit{Ann. Sci.
\'{E}cole Norm. Sup. 6 (1973), 133--252.}

\bigskip

[4] C. D'ANTONI and L. ZSIDO, Groups of linear isometries on multiplier
C*-algebras, \textit{Pacific J. Math. 193 (2000), 279--306.}

\bigskip

[5] J. DIXMIER, Les C*-alg\`{e}bres et leurs repr\'{e}sentations,
Gauthier-Villars, Paris, 1964.

\bigskip

[6] F. FORELLI, Analytic and quasi-invariant measures, \textit{Acta Math. 118
(1967), 33--59.}

\bigskip

[7] F. FORELLI, A maximal algebra, \textit{Math. Scand. 30 (1972), 152--158}.

\bigskip

[8] E. HILLE and R. PHILLIPS, Functional Analysis and Semi-groups, AMS, 1957.

\bigskip

[9] K. HOFFMAN and I. M. SINGER, Maximal subalgebras of $C(\Gamma),$\textit{
Amer. J. Math. 79 (1957), 295-305.}

\bigskip

[10] K. HOFFMAN and I. M. SINGER, Maximal algebras of continuous functions,
\textit{Acta Math. 103 (1960), 217-241.}

\bigskip

[11] A. KISHIMOTO, Simple crossed products of C * -algebras by locally compact
abelian groups, \textit{Yokohama Math. J. 28 (1980), 69--85.}

\bigskip

\bigskip\lbrack12] A. KISHIMOTO, Maximality of the analytic subalgebras of
C*-algebras with flows, \textit{J. Korean Math. Soc. 50 (2013), 1333--1348.}

\bigskip

[13] M. MCASEY, P. MUHLY\ and K.-S. SAITO, Nonselfadjoint crossed products
(invariant subspaces and maximality), \textit{Trans. Amer. Math. Soc. 248
(1979), 381--409.}

\bigskip

[14] P. MUHLY, Function algebras and flows, \textit{Acta Sci. Math. (Szeged)
35 (1973), 111--121.}

\bigskip

[15] D. OLESEN, Inner*-automorphisms of simple C*-algebras, \textit{Comm.
Math. Phys. 44 (1975), 175--190.}

\bigskip

[16] D. OLESEN, G. K. PEDERSEN and E. STORMER, Compact abelian groups of
automorphisms of simple C*-algebras, \textit{Invent. Math. 39 (1977), 55--64.
}

\bigskip

[17] G. K. PEDERSEN, C*-algebras and their automorphism groups, Academic
Press, 1979.

\bigskip

[18] C. PELIGRAD and S. RUBINSTEIN, Maximal subalgebras of C*-crossed
products, \textit{Pacific J. Math. 110 (1984), 325-333.}

\bigskip

[19] C. PELIGRAD and L. ZSIDO, Maximal subalgebras of C*-algebras associated
with periodic flows, \textit{J.Funct. Anal. 262 (2012), 3626-3637.}

\bigskip

[20] W. RUDIN, Fourier Analysis on groups, Interscience, New York, 1962.

\bigskip

[21] D. SARASON, Algebras of functions on the unit circle, \textit{Bull. Amer.
Math. Soc. 79 (1973), 286--299.}

\bigskip

[22] H. SCHAEFER, Topological vector spaces, Springer Verlag New York
Heidelberg Berlin 1971.

\bigskip

[23] A. B. SIMON, On the maximality of vanishing algebras, \textit{Amer. J.
Math. 81 (1959), 613-616.}

\bigskip

[24] B. SOLEL, Algebras of analytic operators associated with a periodic flow
on a von Neumann algebra, \textit{Canad. J. Math. 37 (1985), 405--429.}

\bigskip

[25] B. SOLEL, Maximality of analytic operator algebras, \textit{Israel J.
Math. 62 (1988), 63--89.}

\bigskip

[26] D. C. TAYLOR, The strict topology for double centralizer algebras,
\textit{Trans. Amer. Math. Soc., 150 (1970), 633--643.}

\bigskip

[27] J. WERMER, On algebras of continuous functions, \textit{Proc. Amer. Math.
Soc. 4, (1953). 866--869.}

\bigskip

[28] L. ZSIDO, Spectral and ergodic properties of the analytic generators,
\textit{J. Approximation Theory 20 (1977), 77--138.}

\bigskip

[29] L. ZSIDO, On spectral subspaces associated to locally compact abelian
groups of operators, \textit{Adv. in Math. 36 (1980), 213--276.}

\bigskip

\bigskip

\bigskip

\bigskip
\end{document}